\documentclass[conference]{IEEEtran}
\usepackage{amsmath,amsfonts,amssymb,amsthm}
\usepackage[vlined,algoruled,titlenumbered,noend]{algorithm2e}
\usepackage{array}
\usepackage{epsfig,subfigure,graphicx}
\usepackage[hidelinks]{hyperref}

\def\argmax{\operatornamewithlimits{arg\max}}
\newcommand{\casemax}{\mathrm{casemax}}
\newcommand{\MarsRover}{\textsc{Mars Rover }}
\newcommand{\casemin}{\mathrm{casemin}}
\newcommand{\UB}{\mathit{UB}}
\newcommand{\LB}{\mathit{LB}}
\newcommand{\IND}{\mathit{Ind}}
\newcommand{\CONS}{\mathit{Cons}}
\newcommand{\Root}{\mathit{Root}}
\newcommand{\Max}{\mathit{Max}}
\newcommand{\sq}{\hspace{-1mm}}
\newcommand{\sqm}{\hspace{-2mm}}

\newcommand{\false}{\mathit{false}}

\newcommand{\InventoryControl}{\textsc{Inventory Control }}
\newcommand{\WaterReservoir}{\textsc{Reservoir Management }}

\newtheorem*{example*}{Example}

\newcommand{\comment}[1]{}

\newenvironment{mydef}[1][Definition]{\begin{trivlist}
\item[\hskip \labelsep {\bfseries #1}]}{\end{trivlist}}

\begin{document}

\title{Comparative Analysis of Discrete and Continuous Action Spaces in Reservoir Management and Inventory Control Problems}

\author{
    \IEEEauthorblockN{Sravani Boddepalli}
    \IEEEauthorblockA{
        M.S. in Mathematics,\\
        M.S. in Data Science\\
        Clarkson University, USA\\
        boddeps@clarkson.edu
    }
    \and
    \IEEEauthorblockN{Prathamesh Kothavale}
    \IEEEauthorblockA{
        M.S. in Data Science\\
        Clarkson University, USA\\
        kothavp@clarkson.edu
    }
}

\maketitle

\begin{abstract}
    This paper presents a comparative analysis of discrete and continuous action spaces within the contexts of reservoir management and inventory control problems. We explore the computational trade-offs between discrete action discretizations and continuous action settings, focusing on their effects on time complexity and space requirements across different horizons. Our analysis includes a detailed evaluation of discretization levels in reservoir management, highlighting that finer discretizations approach the performance of continuous actions but at increased computational costs. For inventory control, we investigate deterministic and stochastic demand scenarios, demonstrating the exponential growth in time and space with increasing discrete actions and inventory items. We also introduce a novel symbolic approach for solving continuous problems in hybrid MDPs (H-MDPs), utilizing a new XADD data structure to manage piecewise symbolic value functions. Our results underscore the challenges of scaling solutions and provide insights into efficient handling of discrete and continuous action spaces in complex decision problems. Future research directions include exploring heuristic search methods and improved approximations for enhancing the practicality of exact solutions.
\end{abstract}

\section{Introduction}
\label{Introduction}

Some genuine stochastic arranging issues, for example, those in Mars Wanderer route, stock control, and water supply the executives, include ceaseless factors in their state and activity portrayals. For instance, in Mars Meanderer errands, a wanderer explores inside a consistent spatial climate while performing logical undertakings; in stock control issues for nonstop assets like oil, a business chooses the amount of every thing to arrange in view of unsure interest, limit imperatives, and reordering costs; and in water supply issues, utilities oversee water levels continuously to keep away from undercurrent while improving power age revenue~\cite{bresina02, Scarf_Karlin58, reservoir}.

Definite answers for multivariate consistent state and activity Half and half Markov Choice Cycles (HMDPs) have been restricted. While direct quadratic Gaussian (LQG) control gives precise answers for straight frameworks with Gaussian noise~\cite{lqgc}, most issues including piecewise elements have seen little improvement, particularly in situations where the progress elements include various persistent factors. For discrete activity HMDPs, ideal piecewise straight worth capabilities have been derived~\cite{feng04}, yet broad classes of HMDPs with direct elements for ceaseless factors stay inexplicable.

This work proposes advances to the cutting edge with novel emblematic powerful programming (SDP) strategies, which broaden the class of HMDPs for which careful arrangements can be inferred. These procedures consider taking care of issues with piecewise straight elements and compensations in consistent activity settings and erratic elements in discrete activity settings. As a key application, the methodology empowers the specific arrangement of complicated stock control issues, denoting the initial time in north of 50 years that such an answer has been determined for consistent state situations in this domain~\cite{Scarf_Karlin58}.

\textsc{Discrete Action} \InventoryControl (\textsc{DAIC}): 
A multi-item ($K$-item) inventory consists of continuous amounts of specific items $x_i$ where $i \in [0,K]$ is the number of items and $x_i \in [0,200]$. The customer demand is a stochastic boolean variable $d$ for low or high demand levels.  The order action $a_j$ takes two values of $(0,200)$ where the first indicates no ordering and the second assumes maximum amount of ordering which is 200. There are linear reorder costs and also a penalty for holding items. The transition and reward functions have to be defined for each continuous item $x_i$ and action $j$.

\vspace{2mm}
\textsc{Continuous Action}  \InventoryControl (\textsc{CAIC}):
In a more general  continuous action HMDP setting to this problem, the inventory can order  any of the $i$ items $a_i \in [0,200]$ considering the stochastic customer demand. 

The transition functions for the continuous state $x_i$ and actions $a_i$ is defined as: 
{
\vspace{-2mm}
\begin{figure*}[h]
\centering
\footnotesize
\begin{align}
x'_i &= 
\begin{cases}
d & : x_i + a_i - 150 \\
\neg d & : x_i + a_i - 50    
\end{cases} \notag \\
P(d' = \mathit{true} | d, \vec{x}, \vec{x'}) &= 
\begin{cases}
d & : 0.7 \\
\neg d & : 0.3    
\end{cases} \label{eq:trans_inv}
\end{align}
\end{figure*}
\vspace{2mm}}

The reward is the sum of $K$ functions $R = \sum_{i=0}^K R_i $ as below:
\begin{figure*}[h]
\centering
\footnotesize
\begin{align}
R &= 
\begin{cases}
\sum_{j} x_j \geq C & : -\infty  \\	
\sum_{j} x'_j \geq C & : -\infty 	\\
\sum_{j} x_j \leq C & : 0 \\  
\sum_{j} x'_j \leq C & : 0 	\\
\end{cases}
+ \sum_{i=0}^K \Bigg( 
\begin{cases}
d \wedge x_i \geq 150 & : 150 - 0.1 \cdot a_i - 0.05 \cdot x_i \\
d \wedge x_i \leq 150 & : x_i - 0.1 \cdot a_i - 0.05 \cdot x_i \\
\neg d \wedge x_i \geq 50 & : 50 - 0.1 \cdot a_i - 0.05 \cdot x_i  \\
\neg d \wedge x_i \leq 50 & : x_i - 0.1 \cdot a_i - 0.05 \cdot x_i  \\
\end{cases} \notag \\
&+ \sum_{i=0}^K \Bigg( 
\begin{cases}
x_i \leq 0 & : -\infty  \\	
x'_i \leq 0 & : -\infty  \\	
x_i \geq 0 & : 0 \\
x'_i \geq 0 & : 0
\end{cases} 
\Bigg) \notag
\end{align}
\label{rew_inv}
\end{figure*}

where $C$ is the total capacity for $K$ items in the inventory. The first and last cases check the safe ranges of the capacity such that the inventory capacity of each item above zero and the sum of total capacity below $C$ is desired.

\begin{figure*}[h]
\centering

\begin{subfigure}
                \centering
                \includegraphics[width=0.68\textwidth]{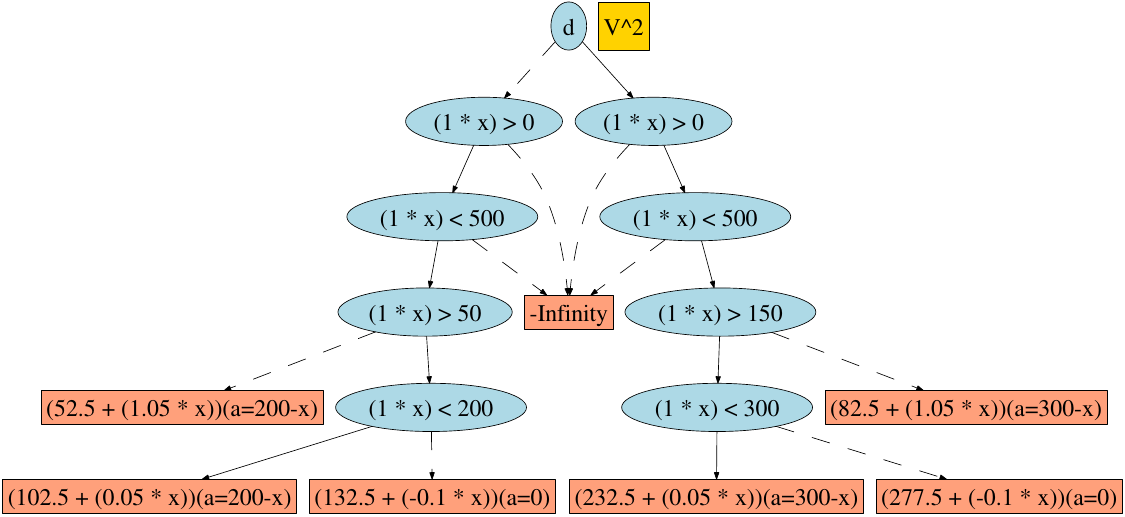}
        \end{subfigure}
                \hspace{2mm}
\begin{subfigure}
                \centering
                \includegraphics[width=0.26\textwidth]{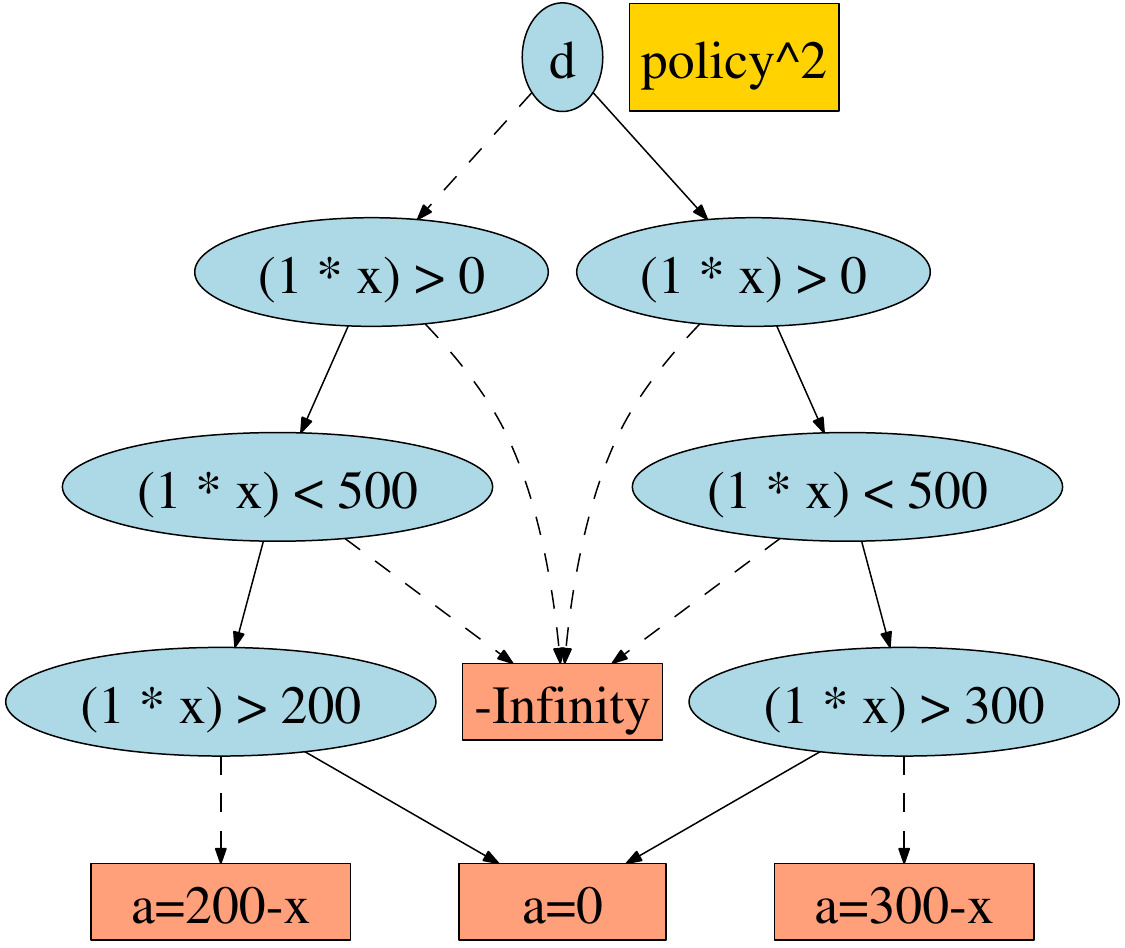}
        \end{subfigure}
\vspace{-2mm}
\caption{\footnotesize The ideal worth capability \( V^2(x) \) for the CAIC issue is addressed by a XADD. To assess \( V^2(x) \), follow the choice tree to a leaf, where the non-coincidental enunciation gives the worth, and the coincidental verbalization gives the ideal strategy \( a = \pi^{*,2}(x) \). The right diagram shows the superior arrangement \( \pi^2 \), lined up with Scarf's answer.}
\label{fig:inv_policy}

\end{figure*}

Note that illegal state values are defined using $-\infty$, in this case having the capacity lower than zero at any time and having capacity higher than that of the total $C$. If our objective is to maximize the long-term \emph{value} $V$ (i.e. the sum of rewards received over an infinite horizon of actions), we show that the optimal value function can be derived in closed-form. 
For a single-item CAIC problem the optimal value function for the second horizon is defined as below:
\vspace{-3mm}
\begin{align}
V = \begin{cases}
(x < 0 \vee x>500) &: -\infty \\
d \land (0 \leq x \leq 500) \land (x \geq 300) &:  277.5 - 0.1 * x \\
d \land (150 \leq x \leq 300) &:  232.5 + 0.05 * x \\
d \land ( 0 \leq x \leq 150) &:  82.5 + 1.05 * x \\
\neg d \land (200 \leq x \leq 500)  &:  132.5 - 0.1 * x \\
\neg d \land (50 \leq x \leq 200) &: 102.5 + 0.05 * x \\
\neg d \land (0 \leq x \leq 50) &:  52.5 + 1.05 * x \\
\end{cases} \label{eq:vfun_inv}
\vspace{-7mm}
\end{align}
The policy obtained from this piecewise and linear value function and $V^2$ itself are shown in Figure~\ref{fig:inv_policy} using an extended algebraic decision diagram (XADD) representation which allows efficient implementation of the \emph{case calculus} for arbitrary functions. According to Scarf's policy for the \InventoryControl problem, if the holding and storage costs are linear the optimal policy in each horizon is always of $(S,s)$ ~\cite{Scarf_Karlin58}. In general this means if ($x>s$) the policy should be not to order any items and if ($x<s$) then ordering $S-s-x$ items is optimal. 

According to this we can rewrite Scarf's policy where each slice of the state space matches with this general rule: 
\begin{align}
\pi^{*,2}(x) = 
\begin{cases}
(x < 0 \vee x>500) &: -\infty \\
d \land (300 \leq x \leq 500)  &:  0 \\
d \land (0 \leq x \leq 300) &:  300 - x \\
\neg d \land (200 \leq x \leq 500) &:  0 \\
\neg d \land (0 \leq x \leq 200) &:  200 - x \\
\end{cases}\nonumber
\end{align}

While this simple example illustrates the power of using continuous variables, for a multi-variate problem it is the very \textit{first solution} to exactly solving problems such as the DAIC and CAIC. 
We propose novel ideas to work around some of the expressiveness limitations of previous approaches, significantly generalizing the range of HMDPs that can be solved exactly.  To achieve this more general solution, this
paper contributes a number of important advances:
\begin{itemize}
\item The use of case calculus allows us to perform Symbolic dynamic programming (SDP) \cite{fomdp} used to solve MDPs with
piecewise transitions and reward functions defined in first-order logic. We define all required operations for SDP such as $\oplus,\ominus,max,min$ as well as new operations such as the continuous maximization of an action parameter $y$ defined as $max_y$ and integration of discrete noisy transition.
\item We perform value iteration for two different settings. In the first setting of DA-HMDP we consider continuous state variables with a discrete action set while in the second setting CA-HMDP we consider continuous states and actions. Both DA-HMDPs and CA-HMDPs are evaluated on various problem domains. The results show that DA-HMDPs applies to a wide range of transition and reward functions providing hyper-rectangular value functions. CA-HMDPs have more restriction in modeling due to the increased complexity caused by continuous actions, and limit solutions to linear and quadratic transitions and rewards but provide strong results for many problems never solved exactly before. 
\item While the \emph{case} portrayal for the ideal \textsc{CAIC}
arrangement displayed in \eqref{eq:vfun_inv} is adequate in principle to
address the ideal worth capabilities that our HMDP arrangement
produces, this portrayal is nonsensical to keep up with by and by
since the quantity of case allotments might become dramatically on
each retreating skyline control step. For \emph{discrete} considered
MDPs, arithmetical choice graphs (ADDs) \cite{bahar93add} have been
effectively utilized in accurate calculations like SPUDD \cite{spudd} to
keep up with reduced esteem portrayals. Roused by this work we
present expanded ADDs (XADDs) to address general minimalistically
piecewise works and tell the best way to perform effective procedure on
them \emph{including} representative expansion. Additionally we present all properties and calculations expected for XADDs.
\end{itemize}

Supported by these algorithmic and information structure progresses, we experimentally show that our SDP approach with XADDs can precisely tackle various HMDPs with discrete and constant activities.

\section{Hybrid MDPs (HMDPs)}
The mathematical framework of Markov Decision Processes (MDPs) is used for modelling many stochastic sequential decision making problems ~\cite{bellman}. This discrete-time stochastic control process chooses an action $a$ available at state $s$. The process then transitions to the next state $s'$ according to $T(s,s')$ and receives a reward $R(s,a)$. The transition function follows the Markov property allowing each state to only depend on its previous state.  We provide novel exact solutions using the MDP framework for discrete and continuous variables in the state and action space. Hybrid state and action MDPs (HMDPs) are introduced in the next section followed by the finite-horizon solution via dynamic programming ~\cite{li05}.
\subsection{Factored Representation}
\label{sec:HMDPs}
In an HMDP, states are represented by vectors of variables \( (\vec{b}, \vec{x}) = (b_1, \ldots, b_n, x_1, \ldots, x_m) \). We assume that each \( b_i \in \{0,1\} \) (\( 1 \leq i \leq n \)) is boolean and each \( x_j \in \mathbb{R} \) (\( 1 \leq j \leq m \)) is continuous. We also assume a finite set of \( p \) actions 
\( A = \{a_1(\vec{y}_1), \ldots, a_p(\vec{y}_p)\} \), where each action \( a_k(\vec{y}_k) \) (\( 1 \leq k \leq p \)) with parameter \( \vec{y}_k \in \mathbb{R}^{|\vec{y}_k|} \) denotes continuous parameters for action \( a_k \), and if \( |\vec{y}_k| = 0 \), then action \( a_k \) has no parameters and is a discrete action.

Each HMDP model requires the following definitions: 

\begin{itemize}
\item  State progress model $P(\vec{b}',\vec{x}'|\vec{b},\vec{x},a,\vec{y})$, which determines the
likelihood of the following state $(\vec{b}',\vec{x}')$ molded on a
subset of the past and next state and activity $a$ with its potential boundaries $\vec{y}$;

\item  Award capability $R(\vec{b},\vec{x},\vec{b}',\vec{x}',a,\vec{y})$, which determines the prompt prize got by making a move $a(\vec{y})$ in state $(\vec{b},\vec{x})$;

\item  Rebate factor $\gamma, \; 0 \leq \gamma \leq 1$ to decide the loads of remunerations in each time step.
\footnote{If time is expressly included as one of the
ceaseless state factors, $\gamma = 1$ is commonly utilized, except if
limiting by skyline (not the same as the state variable time) is
still intended.}
\end{itemize}
A strategy $\pi$ determines the move $a(\vec{y}) =\pi(\vec{b},\vec{x})$ to make in each state $(\vec{b},\vec{x})$. Our
objective is to track down an ideal arrangement of limited skyline subordinate
approaches
\footnote{We expect a limited skyline $H$ in this
paper, but in situations where our SDP calculation joins
in limited time, the subsequent worth capability and
comparing strategy are ideal for $H=\infty$.
For limitedly limited esteem
with $\gamma = 1$, the impending SDP calculation might end in
limited time, however isn't ensured to do as such; for $\gamma < 1$, an
$\epsilon$-ideal strategy for erratic $\epsilon$ can be registered by
SDP in limited time.
}
$\Pi^* = (\pi^{*,1},\ldots,\pi^{*,H})$ that
expands the normal amount of limited compensations more than a skyline $h \in H; H \geq 0$:
\begin{align}
V^{\Pi^*}(\vec{x}) and = E_{\Pi^*} \left[ \sum_{h=0}^{H} \gamma^h \cdot r^h \Big| \vec{b}_0,\vec{x}_0 \right]. \label{eq:vfun_def}
\end{align}
Here $r^h$ is the award gotten at skyline $h$ following $\Pi^*$ where
we expect beginning state $(\vec{b}_0,\vec{x}_0)$ at $h=0$.
 
\begin{figure*}[h]

 \begin{subfigure}
                \centering
                \includegraphics[width=0.2\textwidth]{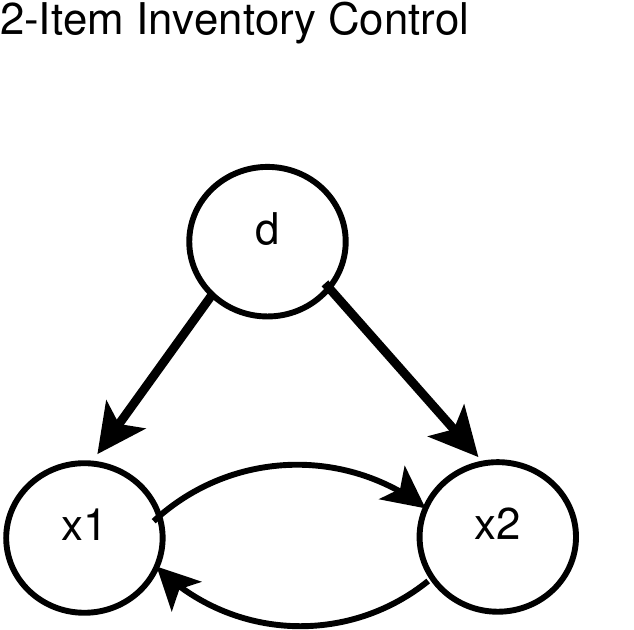}
        \end{subfigure}
                \hspace{2mm}
\begin{subfigure}
                \centering
                \includegraphics[width=0.25\textwidth]{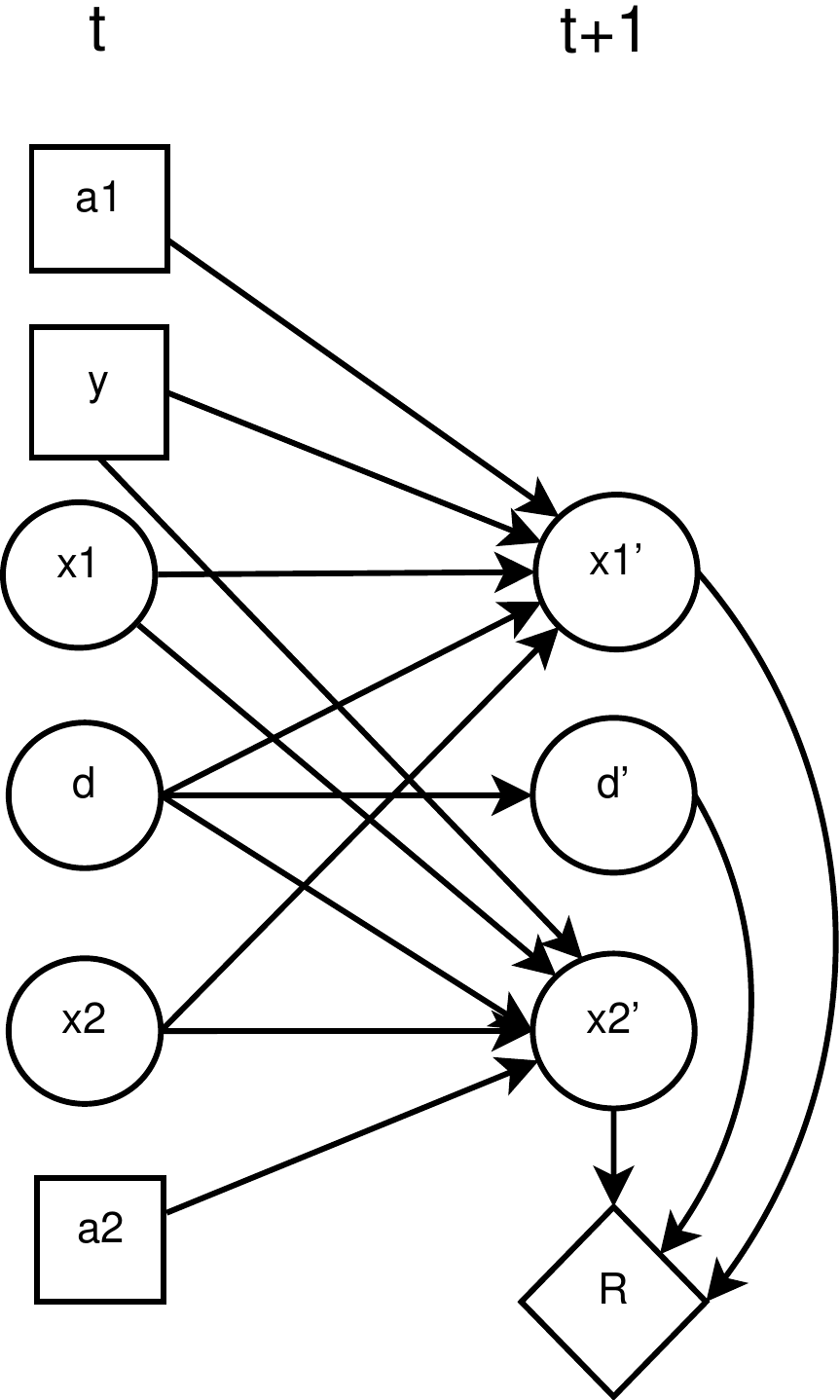}
        \end{subfigure}
                        \hspace{-1mm}
\begin{subfigure}
                \centering
                \includegraphics[width=0.55\textwidth]{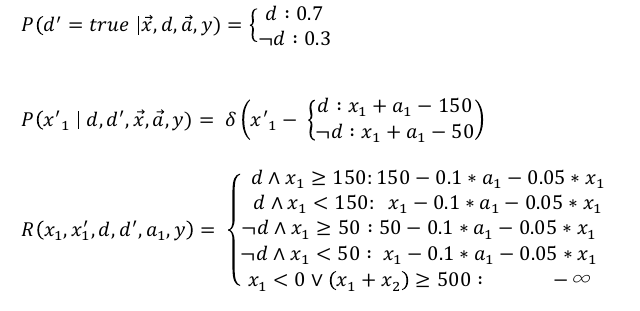}
        \end{subfigure}
        \vspace{-3mm}
\caption{\footnotesize Network geography between state factors in the 2-thing consistent activity \InventoryControl (CAIC) issue (Left); Unique bayes organization (DBN) structure addressing the change and award capability (Center); progress probabilities and prize capability as far as CPF and PLE for $x_1$ (Right).}
\label{fig:dbn}
\vspace{-3mm}
\end{figure*}

Such HMDPs are normally figured \cite{boutilier99dt}
as far as state factors $(\vec{b},\vec{x},\vec{y})$ where possibly $\vec{y} =0$. The change design can be taken advantage of as a unique Bayes
net (DBN) \cite{dbn} where the contingent probabilities
$P(b_i'|\cdots)$ and $P(x_j'|\cdots)$ for each next state variable can
condition on the activity, current and next state.
We can likewise have \emph{synchronic arcs} (factors that condition on each
other in a similar time cut) inside the double $\vec{b}$ or
persistent factors $\vec{x}$ and from $\vec{b}$ to $\vec{x}$.
Thus we can factorize the joint progress model as
{
\begin{align}
P(\vec{b}', \vec{x}' | \vec{b}, \vec{x}, a, \vec{y}) &= \prod_{i=1}^n P(b_i' | \vec{b}, \vec{x}, \vec{b'}, \vec{x'}, a, \vec{y}) \nonumber \\
&\quad \times \prod_{j=1}^m P(x_j' | \vec{b}, \vec{b'}, \vec{x}, \vec{x'}, a, \vec{y}). \nonumber
\end{align}}
where $P(b_i'|\vec{b},\vec{x},\vec{b'},\vec{x'},a,\vec{y})$ might condition on a subset of
$\vec{b}$ and $\vec{x}$ in the current and next state and moreover
$P(x_j'|\vec{b},\vec{b}',\vec{x},\vec{x'},a,\vec{y})$ might condition on a subset of
$\vec{b}$, $\vec{b}'$, $\vec{x}$ and $\vec{x'}$. Figure~\ref{fig:dbn} presents the DBN for a 2-thing CAIC model as per this definition.

We call the restrictive probabilities
$P(b_i'|\vec{b},\vec{x},\vec{b'},\vec{x'},a,\vec{y})$ for \emph{binary} factors $b_i$
($1 \leq I \leq n$) contingent likelihood capabilities (CPFs) - - - not
even lists - - - in light of the fact that overall these capabilities would be able
condition on both discrete and nonstop state as
in the right-hand side of~\eqref{eq:trans_inv}. For the \emph{continuous} factors
$x_j$ ($1 \leq j \leq m$), we address the CPFs
$P(x_j'|\vec{b},\vec{b'},\vec{x},\vec{x'},a,\vec{y})$ with \emph{piecewise
straight equations} (PLEs) fulfilling the accompanying properties:
\begin{itemize}
\item PLEs can condition on the activity, present status, and past state factors
\item  PLEs are deterministic truly intending that to be addressed by probabilities they
should be encoded utilizing Dirac $\delta[\cdot]$ capabilities (model approaching)
\item PLEs are piecewise direct, where the piecewise conditions might be erratic coherent blends of $\vec{b}$, $\vec{b}'$
what's more, straight disparities more than $\vec{x}$ and $\vec{x'}$.
\end{itemize}
The progress capability model gave in the left-hand side of~\eqref{eq:trans_inv} can be communicated in PLE arrangement like the right figure in Figure~\ref{fig:dbn}.
The utilization of the $\delta[\cdot]$ capability guarantees that the PLEs are restrictive
likelihood works that coordinates to 1 more than $x'_j$; In more natural
terms, one can see that this $\delta[\cdot]$ is a straightforward method for encoding
the PLE change $x' = \left\{ \ldots \right.$ as
$P(x_j'|\vec{b},\vec{b'},\vec{x},\vec{x'},a,\vec{y})$.

While it will be certain that our limitations don't allow general stochastic progress commotion (e.g., Gaussian clamor as in LQG control), they truly do allow discrete clamor as in
$P(x_j'|\vec{b},\vec{b'},\vec{x},\vec{x'},a,\vec{y})$ might condition on
$\vec{b'}$, which are stochastically inspected by their CPFs.
\footnote{Continuous stochastic commotion for the progress capability is an on going work which permits us to show stochasticity more generally}
We note that this portrayal successfully permits displaying of
nonstop factor changes as a combination of $\delta$ capabilities,
which has been utilized oftentimes in past precise persistent state MDP
arrangements \cite{feng04,hao09}.
Moreover, we note that our
DA-HMDPs portrayal is more broad than \cite{feng04,li05,hao09} in that
we don't confine the conditions to be straight, but instead
permit it to indicate \emph{arbitrary} capabilities (e.g., nonlinear).

The award capability in DA-HMDPs is characterized as \emph{arbitrary} capability of the present status for each activity $a \in A$.
While experimental models all through the paper will exhibit the full expressiveness of our representative unique programming approach, we note that there are computational benefits to be had when the prize and progress case conditions and works can be confined to straight polynomials.

Because of similar limitations, for CA-HMDPs the prize capability $R(\vec{b},\vec{b'},\vec{x},\vec{x'}, a,\vec{y})$ is characterized as both of the accompanying:

(I) a general piecewise direct capability (boolean or straight circumstances and direct qualities) as in equation~\eqref{rew_inv}; or

(ii) a piecewise quadratic capability of univariate state and a direct capability of univariate activity boundaries:
\begin{align}
R(x,x',d,d', a) and = \begin{cases}
\neg d \land x \geq - 2 \land x \leq 2 : and 4 - x^2 \\
d \lor x < - 2 \lor x > 2 : and 0
\end{cases} \nonumber
\end{align}

These progress and prize requirements will guarantee that all determined capabilities in the arrangement of HMDPs stick to the award
imperatives.

\subsection{Solution methods}
\label{sec:soln}
Presently we give a consistent state speculation of {\it esteem
iteration} \cite{bellman}, which is a powerful programming calculation
for building ideal arrangements. It continues by building a
series of $h$-stage-to-go worth capabilities $V^h(\vec{b},\vec{x})$.
Introducing $V^0(\vec{b},\vec{x}) = 0$ we characterize the quality
$Q_a^{h}(\vec{b},\vec{x},\vec{y})$ of making a move $a(\vec{y})$ in state
$(\vec{b},\vec{x})$ and acting to get
$V^{h-1}(\vec{b},\vec{x})$ from that point as the accompanying:
\vspace{-4mm}
{
\begin{align}
Q_a^{h}(\vec{b},\vec{x},\vec{y}) &= \sum_{\vec{b}'} \int \left( \prod_{i=1}^n P(b_i'|\vec{b},\vec{x},\vec{b'},\vec{x'},a,\vec{y}) \right. \nonumber \\
& \quad \left. \prod_{j=1}^m P(x_j'|\vec{b},\vec{b}',\vec{x},\vec{x'},a,\vec{y}) \right) \nonumber \\
& \times \left[ R(\vec{b},\vec{b'},\vec{x},\vec{x'},a,\vec{y}) + \gamma \cdot V^{h-1}(\vec{b}',\vec{x}') d\vec{x}' \right] \label{eq:qfun}
\end{align}
}
Given $Q_a^h(\vec{b},\vec{x},\vec{y})$ for each $a \in A$ where $\vec{y}$ can likewise be unfilled , we can continue
to characterize the $h$-stage-to-go worth capability as follows:
\begin{align}
V^{h}(\vec{b},\vec{x}) & = \max_{a \in A} \max_{\vec{y} \in \mathbb{R}^{|\vec{y}|}} \left\{ Q^{h}_a(\vec{b},\vec{x},\vec{y}) \right\} \label{eq:vfun}
\end{align}

For discrete activities, augmentation more than the ceaseless boundary $\vec{y}$ is precluded. The $max_{\vec{y}}$ administrator characterized in the following area is expected to sum up arrangements from DA-HMDPs to CA-HMDPs.
On the off chance that the skyline $H$ is limited, the ideal worth capability is
acquired by processing $V^H(\vec{b},\vec{x})$ and the ideal
skyline subordinate strategy $\pi^{*,h}$ at each stage $h$ can be without any problem
decided through $\pi^{*,h}(\vec{b},\vec{x}) = \argmax_a
\argmax_{\vec{y}} Q^h_a(\vec{b},\vec{x},\vec{y})$. If the skyline $H
= \infty$ and the ideal strategy has limitedly limited esteem, then, at that point,
esteem emphasis can end at skyline $h$ if $V^{h} = V^{h-1}$;
then, at that point, $V^\infty = V^h$ and $\pi^{*,\infty} = \pi^{*,h}$.

In DA-HMDPs, we can continuously figure the worth capability in even structure;
notwithstanding, how to process this for HMDPs with remuneration and progress
capability as recently characterized is the target of the representative
dynamic programming calculation that we characterize in the following segment.

\section{Symbolic Dynamic Programming} \label{SDP}

As it's name recommends, emblematic powerful programming (SDP) \cite{fomdp}
is just the most common way of performing dynamic programming (for this situation
esteem emphasis) by means of representative control. While SDP as characterized
in \cite{fomdp} was already just utilized with piecewise
steady capabilities, we currently sum up the portrayal to work with
general piecewise capabilities for HMDPs in this article. Utilizing the \emph{mathematical} meanings of the past segment, we
tell the best way to \emph{compute} equations~\eqref{eq:qfun} and \eqref{eq:vfun}
emblematically.

Before we characterize our answer, in any case, we should officially characterize our
case portrayal and representative case administrators.

\subsection{Case Portrayal and Operations}

All through this article, we will accept that every single representative capability
can be addressed in a \emph{case} structure as follows:
{
\begin{align*}
f = 
\begin{cases}
  \phi_1: & f_1 \\ 
 \vdots&\vdots\\ 
  \phi_k: & f_k \\ 
\end{cases}
\end{align*}
}
Here the $\phi_i$ are sensible formulae characterized over the state
$(\vec{b},\vec{x})$ that can incorporate inconsistent consistent ($\land,\lor,\neg$)
mixes of (a) boolean factors in $\vec{b}$ and (b)
imbalances ($\geq,>,\leq,<$), uniformities ($=$), or disequalities ($\neq$)
where the left and right operands can be \emph{any} capability of at least one
factors in $\vec{x}$.
Each $\phi_i$ will be disjoint from the other $\phi_j$ ($j \neq i$);
anyway the $\phi_i$ may not comprehensively cover the state space, consequently
$f$ may just be a \emph{partial function} and might be indistinct for some
state tasks. Overall we require $f$ to be consistent (counting no discontinuities at segment limits) and all activities save this property. The primary activities expected to perform SDP are given in the accompanying to the case math.

\subsubsection*{Scalar multiplication and Negation}
\emph{Unary operations}, for example, scalar augmentation $c\cdot f$ (for
some consistent $c \in \mathbb{R}$) or refutation $-f$ on case explanations
$f$ are introduced underneath; the unary activity is just applied to each
$f_i$ ($1 \leq I \leq k$).
{
\begin{center}
\begin{tabular}{r c c l}
&

  $c \cdot f = \begin{cases}
    \phi_1  : & c \cdot f_1 \\ 
   \vdots&\vdots\\ 
    \phi_k : & c \cdot f_k \\ 
  \end{cases}$
 &
\vspace{10mm}
  $-f = \begin{cases}
    \neg \phi_1 : & f_1 \\ 
   \vdots&\vdots\\ 
    \neg \phi_k : & f_k \\ 
  \end{cases}$
\end{tabular}
\end{center}
}
\vspace{-9mm} 
\subsubsection*{Binary operations}
Instinctively, to play out a \emph{binary
 operation} on two case articulations, we essentially take the cross-item
of the consistent allotments of each case articulation and play out the
comparing procedure on the subsequent matched allotments. Letting
each $\phi_i$ and $\psi_j$ mean nonexclusive first-request formulae, we can
play out the ''cross-total'' $\oplus$ and ''cross-item'' $\otimes$ of two (anonymous) cases in the
following way:
{\footnotesize 
\begin{center}
\begin{tabular}{r c c c l}
&
\hspace{-6mm} 
  $\begin{cases}
    \phi_1: & f_1 \\ 
    \phi_2: & f_2 \\ 
  \end{cases}$
$\oplus$
&
\hspace{-4mm}
  $\begin{cases}
    \psi_1: & g_1 \\ 
    \psi_2: & g_2 \\ 
  \end{cases}$
&
\hspace{-2mm} 
$ = $
&
\hspace{-2mm}
  $\begin{cases}
  \phi_1 \wedge \psi_1: & f_1 + g_1 \\ 
  \phi_1 \wedge \psi_2: & f_1 + g_2 \\ 
  \phi_2 \wedge \psi_1: & f_2 + g_1 \\ 
  \phi_2 \wedge \psi_2: & f_2 + g_2 \\ 
  \end{cases}$
\end{tabular}
\vspace{4mm}
\hspace{-2mm}
,
\hspace{-2mm}
\begin{tabular}{r c c c l}
&
\hspace{-6mm} 
  $\begin{cases}
    \phi_1: & f_1 \\ 
    \phi_2: & f_2 \\ 
  \end{cases}$
$\otimes$
&
\hspace{-4mm}
  $\begin{cases}
    \psi_1: & g_1 \\ 
    \psi_2: & g_2 \\ 
  \end{cases}$
&
\hspace{-2mm} 
$ = $
&
\hspace{-2mm}
  $\begin{cases}
  \phi_1 \wedge \psi_1: & f_1 * g_1 \\ 
  \phi_1 \wedge \psi_2: & f_1 * g_2 \\ 
  \phi_2 \wedge \psi_1: & f_2 * g_1 \\ 
  \phi_2 \wedge \psi_2: & f_2 * g_2 \\ 
  \end{cases}$
\end{tabular}
\end{center}
}
\normalsize
Similarly, we can perform $\ominus$ by taking away parcel values to acquire the outcome. A few parcels coming about because of
the use of the $\oplus$, $\ominus$, and $\otimes$ administrators
might be conflicting (infeasible); we may basically dispose of such
parts as they are unimportant to the capability esteem.

\subsubsection*{Symbolic maximization}
For SDP, we additionally need to perform amplification over activities for~\eqref{eq:vfun} which is genuinely clear
to characterize:

{
\begin{align}
& \hspace{-4mm} \casemax \Bigg( 
  \begin{cases} 
    \phi_1: & f_1 \\
    \phi_2: & f_2 
  \end{cases} \Bigg), \nonumber \\
& \hspace{-4mm} \begin{cases} 
    \psi_1: & g_1 \\
    \psi_2: & g_2 
  \end{cases} \Bigg) \nonumber \\
& = \begin{cases} 
  \phi_1 \wedge \psi_1 \wedge f_1 > g_1  : & f_1 \\
  \phi_1 \wedge \psi_1 \wedge f_1 \leq g_1 : & g_1 \\
  \phi_1 \wedge \psi_2 \wedge f_1 > g_2  : & f_1 \\
  \phi_1 \wedge \psi_2 \wedge f_1 \leq g_2 : & g_2 \\
  \phi_2 \wedge \psi_1 \wedge f_2 > g_1  : & f_2 \\
  \phi_2 \wedge \psi_1 \wedge f_2 \leq g_1 : & g_1 \\
  \phi_2 \wedge \psi_2 \wedge f_2 > g_2  : & f_2 \\
  \phi_2 \wedge \psi_2 \wedge f_2 \leq g_2 : & g_2 
  \end{cases} \label{symbolicMax}
\end{align}
}
One can confirm that the subsequent case proclamation is still
inside the case language characterized already. From the get go
look this might appear as though a cheat and little is acquired
by this representative skillful deception. Be that as it may, essentially
having a case parcel portrayal that is shut under
boost will work with the shut structure relapse
step that we want for SDP.\ \ Besides, the
XADD that we acquaint later will be capable with exploit the
interior choice design of this
expansion to minimally address it considerably more.
\subsubsection*{Restriction}
In the following activity of \emph{restriction}
we need to limit a capability $f$ to apply just in cases
that fulfill some recipe $\phi$, which we compose as $f|_{\phi}$.
This should be possible by essentially adding $\phi$ to each case parcel
as follows:

{
\begin{center}
\begin{tabular}{r c c l}
&
\hspace{-6mm} 
  $f = \begin{cases}
    \phi_1: & f_1 \\ 
   \vdots&\vdots\\ 
    \phi_k: & f_k \\ 
  \end{cases}$
&

&
\hspace{-2mm}
  $f|_{\phi} = \begin{cases}
    \phi_1 \land \phi : & f_1 \\ 
   \vdots&\vdots\\ 
    \phi_k \land \phi : & f_k \\ 
  \end{cases}$
\end{tabular}
\end{center}
}
Obviously $f|_{\phi}$ possibly applies when $\phi$ holds and is
unclear in any case, subsequently $f|_{\phi}$ is a fractional capability
except if $\phi \equiv \top$.

\subsubsection*{Substitution}
\emph{Symbolic substitution} just takes
a set $\sigma$ of factors and their replacements, e.g.,
$\sigma = \{ x_1'/(x_1 \hspace{-.8mm} + \hspace{-.8mm} x_2), x_2'/x_1^2 \hspace{-.3mm} \exp(x_2) \}$ where
the LHS of $/$ addresses the replacement variable and the
RHS of $/$ addresses
the articulation that ought to be subbed in its place. No factor
happening in any RHS articulation of $\sigma$ can likewise happen in any
LHS articulation of $\sigma$.
We compose the replacement of a non-case capability $f_i$ with $\sigma$
as $f_i\sigma$; for instance, for the $\sigma$ characterized already and
$f_i = x_1' + x_2'$ then $f_i\sigma = x_1 + x_2 + x_1^2 \exp(x_2)$ as
would be normal. We can likewise substitute into case parcels $\phi_j$
by applying $\sigma$ to every imbalance operand; for instance, if
$\phi_j \equiv x_1' \leq \exp(x_2')$ then
$\phi_j \sigma \equiv x_1 + x_2 \leq \exp(x_1^2 \exp(x_2))$.
Having now characterized replacement of $\sigma$ for non-case capabilities $f_i$ and case
parcels $\phi_j$ we can characterize it for case explanations overall:
{
\begin{center}
\begin{tabular}{r c c l}
&
\hspace{-6mm} 
  $f = \begin{cases}
    \phi_1: & f_1 \\ 
   \vdots&\vdots\\ 
    \phi_k: & f_k \\ 
  \end{cases}$
&

&
\hspace{-2mm}
  $f\sigma = \begin{cases}
    \phi_1\sigma: & f_1\sigma \\ 
   \vdots&\vdots\\ 
    \phi_k\sigma: & f_k\sigma \\ 
  \end{cases}$
\end{tabular}
\end{center}
}
\normalsize

One property of replacement is that
on the off chance that $f$ has totally unrelated segments $\phi_i$ ($1 \leq I \leq k$)
then, at that point, $f\sigma$ should likewise have totally unrelated parcels - - -
this understands from the intelligent outcome that
on the off chance that $\phi_1 \land \phi_2 \models \bot$
then $\phi_1\sigma \land \phi_2\sigma \models \bot$.

\subsubsection*{Continuous Integration of the $\delta$-function}
\emph{Continuous Integration} assesses the necessary
minimization $\int_{\vec{x}}$ over the nonstop factors
in a capability $f$. One of the \emph{key novel experiences of SDP} with regards to
HMDPs is that the reconciliation
$\int_{x_j'} \delta[x_j' - g(\vec{x})] f dx_j'$
just \emph{triggers the substitution} $\sigma = \{ x_j'/g(\vec{x}) \}$
on $f$, that is

\begin{align}
\int_{x} \delta[x - g(\vec{x})] f dx' \; = \; f \{x / g(\vec{x}) \} . \label{eq:gen_int}
\end{align}

To perform~\eqref{eq:gen_int} on a more general
representation, we obtain: 
\begin{align*}
    = \begin{cases}
    \phi_1: & f \{ x = g_1 \} \\ 
   \vdots&\vdots\\ 
    \phi_k: & f \{ x= g_k \}  \\ 
  \end{cases}
\end{align*}

Here we note that in light of the fact that $f$ is \emph{already} a case
proclamation, we can basically supplant the single parcel $\phi_i$ with the
numerous parts of $f\{ x/g_i \}|_{\phi_i}$.
This diminishes the \emph{nested} case proclamation back down to a non-settled case
articulation as in the accompanying model:

\begin{align*}
    \begin{cases}
      \phi_1: & 
        \begin{cases}
          \psi_1: & f_{11} \\ 
          \psi_2: & f_{12}  \\ 
        \end{cases} \\
      \phi_2: & 
        \begin{cases}
          \psi_1: & f_{21} \\ 
          \psi_2: & f_{22}  \\ 
        \end{cases} \\
    \end{cases} & \; = \;
        \begin{cases}
          \phi_1 \land \psi_1: & f_{11} \\ 
          \phi_1 \land \psi_2: & f_{12}  \\ 
          \phi_2 \land \psi_1: & f_{21} \\ 
          \phi_2 \land \psi_2: & f_{22}  \\ 
        \end{cases} 
\end{align*}

\subsubsection*{Continuous Maximization}
\emph{Continuous Maximization} of a variable $y$ is characterized as $g(\vec{b},\vec{x}) := \max_{\vec{y}}
\, f(\vec{b},\vec{x},\vec{y})$ where we significantly note that
\emph{the} expanding $\vec{y}$ is a capability
$g(\vec{b},\vec{x})$, thus requiring \emph{symbolic}
obliged advancement. We can modify $f(\vec{b},\vec{x},y)$ through
the accompanying correspondences:
\footnote{The second line guarantees that all unlawful qualities are planned to $-\infty$}
{
\begin{align}
\max_y f(\vec{b},\vec{x},y) & = 
\max_y \casemax_i \, \phi_i(\vec{b},\vec{x},y) f_i(\vec{b},\vec{x},y) \nonumber \\
& = \casemax_i \, \fbox{$\max_y \phi_i(\vec{b},\vec{x},y) f_i(\vec{b},\vec{x},y)$} \label{eq:casemax_max}
\end{align}
}
Since the
$\phi_i$ are commonly disjoint and comprehensive,
$f(\vec{b},\vec{x},y) = \casemax_i \, \phi_i(\vec{b},\vec{x},y) f_i(\vec{b},\vec{x},y)$.

Then on the grounds that
$\max_y$ and $\casemax_i$ are commutative and might be reordered,
we can process $\max_y$ for \emph{each case parcel
individually}. Hence to finish this segment we really want as it were
tell the best way to register a solitary segment emblematically
$\max_y \phi_i(\vec{b},\vec{x},y): f_i(\vec{b},\vec{x},y)$.

In $\phi_i$, we see that each conjoined limitation serves one of
three purposes:
\begin{itemize}
\item  \emph{upper bound ($\UB$) on $y$}: can be composed as $y < \cdots$ or $y \leq \cdots$
\item  \emph{lower bound ($\LB$) on $y$}: it very well may be composed as $y >\cdots$ or $y \geq \cdots$
\footnote{For reasons for assessing
a case capability $f$ at an upper or lower bound,
it doesn't make any difference whether a bound is comprehensive ($\leq$ or $\geq$)
or then again selective ($<$ or $>$) since $f$ is expected to be persistent
furthermore, consequently assessing at the constraint of the comprehensive bound will
match the assessment for the restrictive bound.}
\item  \emph{independent of $y$ ($\IND$)}: the imperatives don't contain $y$
what's more, can be securely calculated beyond the $\max_y$.
\end{itemize}

Since there are various emblematic upper and lower
limits on $y$, overall we should apply the $\casemax$
($\casemin$) administrator to decide the most noteworthy lower bound $\LB$
(most minimal upper bound $\UB$).

We likewise know that $\max_y \phi_i(\vec{b},\vec{x},y)
f_i(\vec{b},\vec{x},y)$ for a ceaseless capability $f_i$ should happen at the basic marks of the capability - - -
either the upper or lower limits ($\UB$ and $\LB$) of $y$,
or on the other hand the $\Root$ (i.e., zero) of $\frac{\partial}{\partial y} f_i$
w.r.t.\ $y$. Every one of $\UB$, $\LB$, and $\Root$
is an emblematic capability of $\vec{b}$ and $\vec{x}$.

Given the \emph{potential} maxima points of $y = \UB$, $y = \LB$, and
$y = \Root$ of $\frac{\partial}{\partial y} f_i(\vec{b},\vec{x},y)$
w.r.t. requirements $\phi_i(\vec{b},\vec{x},y)$ - - - which are all
emblematic capabilities - - - we should emblematically assess which yields the
expanding esteem $\Max$ for this case parcel:

Here $\casemax(f,g,h) = \casemax(f,\casemax(g,h))$. The
replacement administrator $\{ y/f \}$ replaces $y$ with case proclamation $f$,
characterized already.

Right now, we have nearly finished the calculation
of the $\max_y \phi_i(\vec{b},\vec{x},y) f_i(\vec{b},\vec{x},y)$
with the exception of one issue: the joining of the autonomous ($\IND$) limitations
(figured out beforehand) and extra requirements that emerge from the
representative nature of the $\UB$, $\LB$, and $\Root$.

Explicitly for the last option, we really want to guarantee that for sure $\LB \leq \Root \leq \UB$
(or on the other hand in the event that no root exists, $\LB \leq \UB$) by building a set
of imperatives $\CONS$ that guarantee these circumstances hold; to do this,
it gets the job done to guarantee that for every conceivable articulation $e$ used to
build $\LB$ that $e \leq \Root$ and comparably for the $Root$ and $\UB$.
Presently we express the end-product as a solitary case parcel:
\begin{equation*}
\max_y \phi_i(\vec{b},\vec{x},y) f_i(\vec{b},\vec{x},y) \;\; = \;\;
\left\{ \CONS \land \IND: \Max \right.
\end{equation*}
Thus, to finish the expansion for a whole case proclamation $f$, we really want just apply the above technique to each case segment of $f$ and afterward play out a representative $\casemax$ on the outcomes in general.

\subsection{Symbolic Dynamic Programming (SDP)}

\incmargin{.5em}
\linesnumbered
\begin{algorithm}[h]
\vspace{-.5mm}
\dontprintsemicolon
\SetKwFunction{regress}{Regress}
\Begin
{
   $V^0:=0, h:=0$\;
   \While{$h < H$}
   {
       $h:=h+1$\;
       \ForEach {$a(\vec{y}) \in A$}
       {
              $Q_a^{h}(\vec{y})\,:=\,$\regress{$V^{h-1},a,\vec{y}$}\;
			  \emph{//Continuous action parameter}\;
			  \If  {$ \mid \vec{y} \mid >0$}  
               {
               $Q_a^{h}(\vec{y}) := \max_{\vec{y}} \, Q_a^{h}(\vec{y})$ $\,$ \;
               $\pi^{*,h} := \argmax_{a} \, Q_a^{h}(\vec{y})$\;
               } 
               \Else 
               { $\pi^{*,h} := \argmax_{a} \, Q_a^{h}(\vec{y})$ \; }
        }
       $V^{h} := \casemax \, Q_a^{h}(\vec{y})$ $\,$ \;
       \If{$V^h = V^{h-1}$}
           {break $\,$ \emph{// Terminate if early convergence}\;}
   }
     \Return{$(V^h,\pi^{*,h})$} \;
}
\caption{\footnotesize \texttt{VI}(HMDP, $H$) $\longrightarrow$ $(V^h,\pi^{*,h})$ \label{alg:vi}}
\vspace{-1mm}
\end{algorithm}
\decmargin{.5em}

In this section the symbolic value iteration algorithm (SVI) for HMDPs is presented.
Our objective is to take a DA-HMDP or CA-HMDP as defined in Section~\ref{sec:HMDPs}, apply value
iteration as defined in Section~\ref{sec:soln}, and produce
the final value optimal function $V^h$ at horizon $h$ in the form
of a case statement presented in Algorithm~\ref{alg:vi}. 
We use the \textsc{CAIC} example from the introduction to help clarify each step of this algorithm. 

For the base case of $h=0$ in line 2, we note that setting $V^0(\vec{b},\vec{x}) = 0$
(or to the reward case statement, if it is not action dependent)
is trivially in the form of a case statement.

\incmargin{.5em}
\linesnumbered
\begin{algorithm}[h]
\vspace{-.5mm}
\dontprintsemicolon
\SetKwFunction{remapWithPrimes}{Prime}

\Begin{
    $Q=$ \remapWithPrimes{$V$} $\,$ \emph{// All $v_i \to v_i'$ ($\equiv$ all $b_i \to b_i'$ and all $ x_i \to x_i'$)} \;

\If {$v'$ in $R$}
	 {$Q := R(\vec{b},\vec{b}',\vec{x},\vec{x}',a,\vec{y}) \oplus (\gamma \cdot Q)$} \;
    \ForEach { $v'$ in $Q$}  
    {
    	\If {$v'$ = $x'_j$}
    	{
    	\emph{//Continuous marginal integration}\\
         $Q := \int Q \otimes P(x_j'|\vec{b},\vec{b}',\vec{x},\vec{x}',a,\vec{y}) \, d_{x'_j}$\;
    	}
	    \If {$v'$=$b'_i$}
    	{
    	\emph{// Discrete marginal summation}\\
         $Q := \left[ Q \otimes P(b_i'|\vec{b},\vec{b}',\vec{x},\vec{x}',a,\vec{y}) \right]|_{b_i' = 1}$
         $\oplus \left[ Q \otimes P(b_i'|\vec{b},\vec{b},\vec{x},\vec{x}',a,\vec{y}) \right]|_{b_i' = 0}$\;
    	}
    }
    \If {$\neg$ ($v'$  in $R$)}
    {$Q := R(\vec{b},\vec{b}',\vec{x},\vec{x}',a,\vec{y}) \oplus (\gamma \cdot Q)$ }\;
    \Return{$Q$} \;
}
\caption{\footnotesize \texttt{Regress}($V,a,\vec{y}$) $\longrightarrow$ $Q$ \label{alg:regress}}
\vspace{-1mm}
\end{algorithm}
\decmargin{.5em}

Then, for $h > 0$ and for each activity in line 5 we should perform lines 6- - 12. Beginning with the utilization of Algorithm~\ref{alg:regress}. Note that we have discarded boundaries $\vec{b}$ and
$\vec{x}$ from $V$ and $Q$ to stay away from notational mess.
Luckily, given our recently characterized
activities, SDP is direct and can be separated into five
steps:
\begin{enumerate}
\item  {\it Prepare Function}: Since $V^{h}$ will turn into
the ''following state'' in esteem emphasis, we arrangement a replacement
$\sigma = \{ b_1/b_1', \ldots, b_n/b_n', x_1/x_1', \ldots, x_m/x_m' \}$
what's more, get $V'^{h} = V^{h}\sigma$ in line 2 of Algorithm~\ref{alg:regress}. Beginning with the principal emphasis, for the \textsc{CAIC} model this step doesn't make a difference to $h-1=0$ since $V^0=0$.

\item  {\it Add Award Function}: Assuming the prize capability $R$ contains any prepared state variable $b'$ or $x'$, lines 3- - 4 of Algorithm~\ref{alg:regress} is executed to add this award capability to the past limited Q-esteem. In the event that $R$ had no prepared factors, it is added to the Q-esteem toward the finish of Algorithm~\ref{alg:regress} in lines 14- - 15. The award capability of \textsc{CAIC} contains prepared $x'$ subsequently the Q-esteem is characterized as underneath:
\vspace{-1mm}
\begin{figure*}[h]
    \centering
    \begin{align}
    Q = 
    \begin{cases}
    (x < 0 \vee x > 500 \vee x' < 0 \vee x' > 500) & : -\infty \\
    d \land (150 \leq x \leq 500) & : 150 - 0.1 \cdot a - 0.05 \cdot x \\
    d \land (0 \leq x \leq 150) & : 0.95 \cdot x - 0.1 \cdot a \\
    \neg d \land (50 \leq x \leq 500) & : 50 - 0.1 \cdot a - 0.05 \cdot x \\
    \neg d \land (0 \leq x \leq 50) & : 0.95 \cdot x - 0.1 \cdot a
    \end{cases}
    \end{align}
    \vspace{-4mm}
    \caption{Piecewise function $Q$ as a function of $x$, $x'$, and $a$ under different conditions.}
    \label{fig:piecewise_function}
\end{figure*}

\item  {\it Nonstop Integration}: As characterized in line 7- - 9 of Algorithm~\ref{alg:regress} once we have our prepared worth
capability $V'^{h}$ on the off chance that explanation design characterized over next state
factors $(\vec{b}',\vec{x}')$, we assess the indispensable
minimization $\int_{\vec{x}'}$ over the constant factors
in~\eqref{eq:qfun}. Since the lower and upper coordination limits
are separately $-\infty$ and $\infty$
what's more, we have prohibited synchronic curves between factors in $\vec{x}'$
in the progress DBN, we can minimize out each
$x_j'$ freely, and in any request. Concurring to~\eqref{eq:gen_int} we have the accompanying:
\begin{align*}
\int_{x_j'} \delta[x_j' - g(\vec{x})] V'^{h} dx_j' \; = \; V'^{h} \{x_j' / g(\vec{x}) \}  
\end{align*}
This
activity is performed over and over in succession \emph{for each}
$x_j'$ ($1 \leq j \leq m$) for each activity $a$. The as it were
unexpected confusion is that the type of
$P(x_j'|\vec{b},\vec{x},\vec{b'},\vec{x'},a,\vec{y})$ is a \emph{conditional}
condition like the right-hand of Figure~\ref{fig:dbn}, and addressed conventionally
as follows:
\begin{align}
   P(x_j'|\vec{b},\vec{x},\vec{b'},\vec{x'},a,\vec{y}) = \delta\left[ x_j' = \begin{cases}
    \phi_1: & f_1 \\ 
   \vdots&\vdots\\ 
    \phi_k: & f_k \\ 
  \end{cases} \right] \label{eq:cond_sub}
\end{align}
Basically, we can read~\eqref{eq:cond_sub} as a \emph{conditional
substitution}, i.e., in each of the different \emph{previous state}
conditions $\phi_i$ ($1 \leq I \leq k$), we get an alternate
replacement for $x_j'$ showing up in $V'^{h}$ (i.e., $\sigma = \{ x_j'/f_i
\}$).

To play out the full nonstop incorporation,
assuming that we instate
$\tilde{Q}_a^{h+1} := V'^{h}$ for each activity $a \in A$, and rehash
the above integrals for all $x_j'$, refreshing $\tilde{Q}_a^{h+1}$ each time,
then after end of all $x_j'$ ($1 \leq j \leq m$), we will have
the halfway relapse of $V'^{h}$ for the nonstop factors for
each activity $a$ indicated by $\tilde{Q}_a^{h+1}$. Following the \textsc{CAIC} model, constant coordination of $x$ brings about the accompanying:
{\footnotesize
\begin{figure*}[h]
    \centering
    \begin{align}
    Q = 
    \begin{cases}
    x < 0 \vee x > 500 & : -\infty \\
    d \land (x \geq 150) \land (150 \leq (x + a) \leq 650) & : 150 - 0.1 \cdot a - 0.05 \cdot x \\
    d \land (x \geq 150) \land ((x + a \geq 650) \vee (x + a \leq 150)) & : -\infty \\
    d \land (x \leq 150) \land (150 \leq (x + a) \leq 650) & : 0.95 \cdot x - 0.1 \cdot a \\
    d \land (x \leq 150) \land ((x + a \geq 650) \vee (x + a \leq 150)) & : -\infty \\
    \neg d \land (x \geq 50) \land (50 \leq (x + a) \leq 550) & : 50 - 0.1 \cdot a - 0.05 \cdot x \\
    \neg d \land (x \geq 50) \land ((x + a \geq 550) \vee (x + a \leq 50)) & : -\infty \\
    \neg d \land (x \leq 50) \land (50 \leq (x + a) \leq 550) & : 0.95 \cdot x - 0.1 \cdot a \\
    \neg d \land (x \leq 50) \land ((x + a \geq 550) \vee (x + a \leq 50)) & : -\infty
    \end{cases}
    \end{align}
    \vspace{-4mm}
    \caption{Piecewise function $Q$ as a function of $x$, $a$, and $d$ under different conditions.}
    \label{fig:piecewise_function_recent}
\end{figure*}
}
\item {\it Discrete Marginalization}: Now that we have our partial
regression $\tilde{Q}_a^{h+1}$ for each action $a$, we proceed
to derive the full backup $Q_a^{h+1}$ from $\tilde{Q}_a^{h+1}$
by evaluating the discrete 
marginalization $\sum_{\vec{b}'}$ in~\eqref{eq:qfun} which is shown in lines 10--12 of Algorithm~\ref{alg:regress}.
Because we previously disallowed synchronic arcs
between the variables in $\vec{b}'$ 
in the transition DBN, we can sum out each variable $b_i'$ ($1 \leq i \leq n$) 
independently.  Hence, initializing
$Q_a^{h+1} := \tilde{Q}_a^{h+1}$
we perform the discrete regression by applying the following iterative
process \emph{for each} $b_i'$ in any order for each action $a$:
\begin{align}
Q_a^{h+1} := & \left[ Q_a^{h+1} \otimes P(b_i'|\vec{b},\vec{x},\vec{b'},\vec{x'},a,\vec{y}) \right]_{b_i' = 1} \\
             & \oplus \left[ Q_a^{h+1} \otimes P(b_i'|\vec{b},\vec{x},\vec{b'},\vec{x'},a,\vec{y}) \right]_{b_i' = 0}.
\end{align}
This requires a variant of the earlier restriction operator $|_v$ that
actually \emph{sets} the variable $v$ to the given value if present.
Note that both $Q_a^{h+1}$ and $P(b_i'|\vec{b},\vec{x},\vec{b'},\vec{x'},a,\vec{y})$ can be represented
as case statements (discrete CPTs \emph{are} case statements), 
and each operation produces a case statement.
Thus, once this process is complete, we have marginalized over
all $\vec{b}'$ and $Q_a^{h+1}$ is the symbolic representation
of the intended Q-function. In \textsc{CAIC} discrete marginalization of the boolean state variable $d$ is not performed since there is no primed version of this variable $d'$ in the current Q-function. 

\incmargin{.5em}
\linesnumbered
\begin{algorithm}[h]
\vspace{-.5mm}
\SetKwFunction{solveForVar}{SolveForVar}
\dontprintsemicolon
\Begin{
	$LB $ =$- \infty$ , 
	$UB$= $+ \infty$, 
	$IND,CONS $=$\mathit{true}$ , 
	$\mathit{Case_{max}} \,=\, \emptyset$ \;
     \For {$\phi_i \in f$ (\emph{For all partitions of $f$})}  
     {
     	\For{ $c_i \in \phi_i$ (\emph{For all conditions $c$ of $\phi_i$})}
     	{     		
     		\If {$c_i \leq y $} {$LB \,:=\, casemax(LB, c_i )$ \emph{//Add  $c_i$ to  LB, take max of all LBs}\; }
     		\If {$c_i \geq y $}	{$UB \,:=\, casemin(UB,c_i )$ \emph{//Add $c_i$ to  UB, take min of all UBs} \;	}
     		\Else {$IND \,:=\, [IND,c_i ]$ \emph{//Add constraint $c_i$ to  independent constraint set}\;}
     	}    	
     }
	$\mathit{Root}$ :=$ \solveForVar(y, f_i)$ \;
     	
     	\If {($\mathit{Root} \neq \mathit{null}$)}
     	 {
     	 		$CONS = (\mathbb{I} \left[ LB \right] \leq \mathbb{I} \left[ \mathit{Root} \right] ) 
     	 		\wedge (\mathbb{I} \left[ \mathit{Root} \right] \leq \mathbb{I} \left[ UB \right] )$\;
     	 } 
     	 \Else
     	  {
     	 	$CONS \,:=\, (\mathbb{I} \left[ LB \right] \leq \mathbb{I} \left[ UB \right])$\;
     	 }
     	\emph{//Conditions and value of continuous max for this partition\; }   
     	$\mathit{Max} \,:=\,  IND \wedge CONS : \casemax (f_i\left \lbrace y/LB \right \rbrace, f_i \left \lbrace y/UB \right \rbrace 		,f_i \left \lbrace y/ \mathit{Root} \right \rbrace)$\;
     	\emph{//Take maximum of this partition and all other partitions\;}       
     	$\mathit{Case_{max}} \,:=\, max(\mathit{Case_{max}} ,\mathit{Max}) $  
        
     \Return{$\mathit{Case_{max}} $} \;
}
\caption{\footnotesize \texttt{Continuous Maximization}($y$, $f(\vec{b},\vec{x},y)$) $\longrightarrow(max_{y}f(\vec{b},\vec{x},y))$ \label{alg:contMax}}
\vspace{-1mm}
\end{algorithm}
\decmargin{.5em}

\item {\it Continuous action Maximization}: This maximization is over an action variable $a(\vec{y})$ in line 8--9 of
Algorithm~\ref{alg:vi} where $\mid \vec{y} \mid>0$, requires a continuous maximization. Here we take the maximum over parameter $y$ of action variable $a(\vec{y})$. If the action is discrete $\mid \vec{y} \mid=0$ , lines 8--10 are not performed.
Exploiting the commutativity of $\max$, we can first
rewrite any multivariate $\max_{\vec{y}}$ as a sequence of univariate
$\max$ operations $\max_{y_1} \cdots \max_{y_{|\vec{y}|}}$; hence it
suffices to provide just the \emph{univariate} $\max_y$ solution:
\begin{align}
\max_{\vec{y}} =\max_{y_1} \cdots \max_{y_{|\vec{y}|}} \Rightarrow g(\vec{b},\vec{x}) := \max_{y} \, f(\vec{b},\vec{x},y). \nonumber
\end{align}
According to the properties on the Continuous Maximization operation defined in the previous section, we compute a univariate maximization  
$\max_y \phi_i(\vec{b},\vec{x},y) f_i(\vec{b},\vec{x},y)$ using Algorithm \ref{alg:contMax}. 
\footnote{Note also that from here out we assume that all case partition conditions $\phi_i$ of
$f$ consist of conjunctions of non-negated linear inequalities and
possibly negated boolean variables --- conditions easy to enforce
since negation inverts inequalities, e.g., $\neg [x < 2] \equiv [x \geq 2]$
and disjunctions can be split across multiple non-disjunctive, 
disjoint case partitions.} 
Each step of this algorithm is followed using one of the partitions of the Q-function in this case the first partition with the constraints of 
$\phi_i(x,d,a) \equiv d \land (0 \leq x \leq 500) \land (x \geq 150) \land ((x+a) \leq 650) \land ((x+a) \geq 150) $ and function value of $ f_i(x,d,a) = 150 - 0.1 * a - 0.05 * x$.

To begin the set of lower bound $\LB$ is set to $-\infty$ and upper bound $\UB$ to $\infty$ so that any value larger than $-\infty$ is defined as the lower bound and any value lower than $+\infty$ is defined for $\UB$. Constraint variables $IND$ and $CONS$ are assumed to be true and the result of the $\casemax$ is set to empty. 

Each constraint $c_i$ in each partition $\phi_i$ is added to one of the sets of lower bound, upper bound or independent constraint as determined in 5--10. In our example this is equal to $\LB = (150 - x, 0) , \UB= (650 - x, 1000000) $ and $\IND=(d,x \geq 0 , x \leq 500) $ where the $0$ and $1000000$ are the natural lower and upper bounds on any inventory item $x$.
A unique $\LB$ and $\UB$ is defined by taking the maximum of the lower bounds and the minimum of the upper bounds as the best bounds in the current partition and the function \emph{SolveForVar} of line 13 takes any roots of the partition function (not applicable in the current partition)

The boundary constraints in lines 14--17 are added to the independent constraints as the constraint of the final maximum $\mathit{Max}$:
{
\begin{align}
\CONS_{\LB \leq \UB} \sq = & \, \sq [0 \leq 1000000] \\
& \land [150 - x \leq 1000000] \\
& \land [150 - x \leq 650 - x] \\
& \land [0 \leq 650 - x].
\end{align}}
Here, two constraints are tautologies and may be removed.
A $\casemax$ is performed on the substituted $\LB$,$\UB$ and the roots on the function $f_i$: 
{\footnotesize 
\begin{figure*}[h]
\centering
\begin{align}
\Max &= \casemax \Big( f_i \{ y / \Root \} = \textit{null} , \nonumber \\
f_i \{ y / \LB \} &= \begin{cases}
x > 150: & \sqm 150 - 0.1 \cdot (0) - 0.05 \cdot x = 150 - 0.05 \cdot x \\ 
x \leq 150: & \sqm 150 - 0.1 \cdot (150 - x ) - 0.05 \cdot x = 135.00075 + 0.05 \cdot x\\
\end{cases} , \nonumber \\
f_i \{ y / \UB \} &= \begin{cases}
x > -1000000: & \sqm 150 - 0.1 \cdot (650 -x) - 0.05 \cdot x = 84.980494 + 0.05 \cdot x \\ 
x \leq -1000000: & \sqm 150 - 0.1 \cdot (1000000) - 0.05 \cdot x = -99850 - 0.05 \cdot x\\
\end{cases} \Big) \nonumber
\end{align}
\caption{Max function calculation with case analysis.}
\label{fig:max_case_analysis}
\end{figure*}}
Taking this $\max_y$ is performed in line 19 for each partition using both independent and boundary constraints where the resulting maximum is according to the $\casemax$ operator defined in the previous section and the partition of $x \leq -1000000:  -99850 -0.05 * x$ is omitted due to inconsistency.
{\footnotesize 
\begin{align*}
\Max = 
\begin{cases}
(x > 150) \land (x \leq 650):    & \sqm 150 - 0.05 * x \\ 
(x > 150) \land (x \geq 650):    & \sqm 84.980494 + 0.05 * x \\ 
x \leq 150: & \sqm 135.00075 + 0.05 * x\\ 
\end{cases}
\end{align*}}
Returning to~\eqref{eq:casemax_max}, we have now specified the inner operation (shown in the $\Box$) which is defined as the following for this partition.\footnote{ These last two results are defined by taking out all inconsistent partitions. This is done using efficient pruning techniques mentioned in the next section.}
\begin{align*}
\Max & = 
\begin{cases}
d \land (150 \leq x \leq 500):    & \sqm 150 - 0.05 * x \\ 
\text{\normalsize otherwise} : & -\infty\\ 
\end{cases}
\end{align*}
To complete the maximization for an entire case statement $f$, we need only apply the above procedure to each case partition of $f$ and then $\casemax$ all of these results in line 21: 
\begin{align*}
Q & = 
\begin{cases}
d \land (150 \leq x \leq 500):    & \sqm 150 - 0.05 * x \\ 
d \land (0 \leq x \leq 150):    & \sqm -14.99925 + 1.05 * x \\ 
\neg d \land (50 \leq x \leq 500):    & \sqm 50 - 0.05 * x \\ 
\neg d \land (0 \leq x \leq 50):    & \sqm -5 + 1.05 * x \\ 
\text{\normalsize otherwise} : & -\infty\\ 
\end{cases}
\end{align*}
To obtain the policy in Figure~\ref{fig:inv_policy}, 
we need only annotate leaf values with any 
$\UB$, $\LB$, and $\Root$ substitutions performing line 10 or 12 in Algorithm~\ref{alg:vi}.
Continuous maximization is further explained in the next section using the appropriate data structure. 

\item {\it Maximization}: Now that we have $Q_a^{h+1}(\vec{y})$ in
case format for each action $a \in\{a_{1}(\vec{y}_1), \ldots \\, a_{p}(\vec{y}_p)\}$, obtaining
$V^{h+1}$ in case format as defined in~\eqref{eq:vfun} requires
sequentially applying
\emph{symbolic maximization} in line 14 as defined previously:
\begin{align}
V^{h+1} &= \max\left( Q_{a_1}^{h+1}(\vec{y}), \max \left( \ldots, \right. \right. \nonumber \\
         & \quad \max\left( Q_{a_{p-1}}^{h+1}(\vec{y}), Q_{a_p}^{h+1}(\vec{y}) \right) \nonumber \\
         & \left. \left. \ldots \right) \right) \nonumber
\end{align}

\end{enumerate}
Note that for our \textsc{CAIC} model the last Q-capability is equivalent to the ideal worth capability since we have thought about a solitary constant activity here.
By enlistment, in light of the fact that $V^0$ is a case explanation and applying
SDP to $V^h$ on the off chance that explanation structure produces $V^{h+1}$ in the event that
proclamation structure, we have accomplished our planned
objective with SDP. On the issue of rightness,
we note that every activity above basically executes one of the
dynamic programming activities in \eqref{eq:qfun} or \eqref{eq:vfun},
so rightness essentially follows from confirming (a) that each case
activity creates the right outcome and that (b) each case activity
is applied in the right succession as characterized in \eqref{eq:qfun} or
\eqref{eq:vfun}.

On a last note, we see that SDP holds for \emph{any} representative
case articulations; we have not confined ourselves to rectangular
piecewise capabilities, piecewise straight capabilities, or even piecewise
polynomial capabilities. As the SDP arrangement is simply emblematic,
SDP applies to \emph{any} HMDP utilizing limited representative capability
that can be written in the event that configuration! Obviously, that is the hypothesis,
next we meet practice.

\section{Extended Algebric Decision Diagrams (XADDs)} \label{XADD}

In the past segment all activities expected to perform SDP calculations were covered. The case proclamations address inconsistent piecewise capabilities permitting general answers for nonstop issues. By and by, it very well may be restrictively costly to keep up with
a case explanation portrayal of a worth capability with express
parcels. Roused by the SPUDD~\cite{spudd} calculation which
keeps up with reduced esteem capability portrayals for limited discrete
calculated MDPs utilizing arithmetical choice graphs (ADDs)~\cite{bahar93add},
we stretch out this formalism to deal with nonstop factors in an information
structure we allude to as the XADD.

\begin{figure}[h]
\centering

\includegraphics[width=0.8\linewidth]{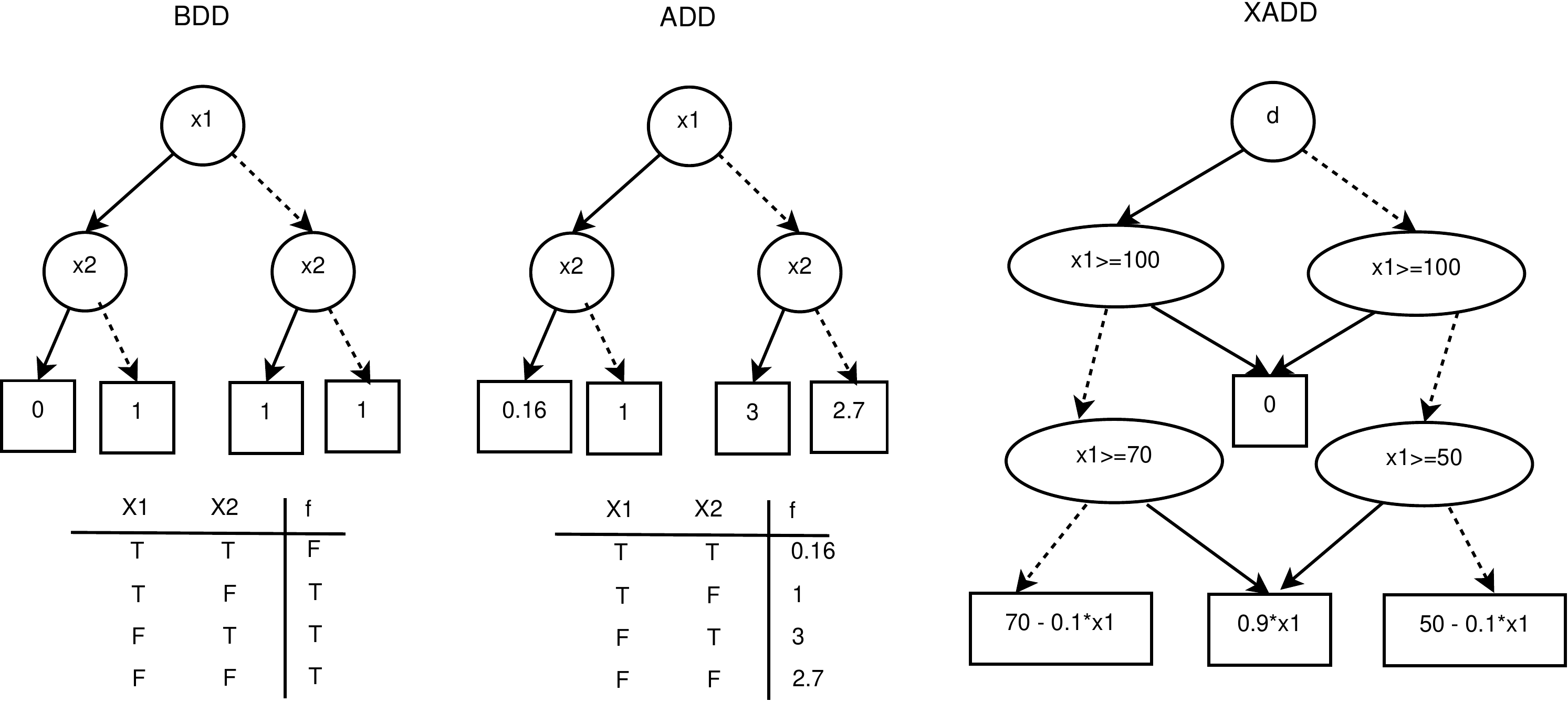}

\caption{\footnotesize Comparison of the three decision diagrams: Binary decision diagrams (BDDs) with boolean leaves and decisions (Left) representing as shown in the truth table; Algebraic decision diagrams (ADDs) with  boolean decision nodes and real values at the leaves (Middle) represented by the truth table;}
\label{fig:bdd_add_xadd}

\end{figure}

Here we present this reduced information design of XADDs which can execute case explanations effectively. Figure~\ref{fig:inv_policy} of the presentation area exhibits the worth capability for the \InventoryControl issue as a XADD portrayal. While XADDs are stretched out from ADDs, ADDs are reached out from Twofold choice charts (BDDs), permitting first-request rationale rather than boolean rationale. Figure~\ref{fig:bdd_add_xadd} exhibits instances of the three choice outlines of BDD, ADD and XADD as a correlation with show their expressiveness.  

A \emph{binary decision diagram} (BDD) ~\cite{bryant} can represent propositional formulas or boolean functions $\lbrace 0,1\rbrace^n \rightarrow \lbrace 0,1\rbrace$ as an ordered \emph{directed acyclic graph} (DAG) where each node represents a random variable and edges represent direct correlations between the variables. Each decision node is a boolean test variable with two successor nodes of false/true. The edge from the decision node to a false (true) child represents assigning  0 (1) in boolean logic. To evaluate the boolean function of a certain BDD, each of the variables are assigned a false/true value by following the corresponding branches until reaching a leaf. The boolean value at the leaf is the value returned by this function according to the given variable assignment. 

Extending BDDs to \emph{algebric decision diagrams} (ADDs) allows a real-value range in the function representation $\lbrace 0,1\rbrace^n \rightarrow \mathbb{R}$. ADDs further provide an efficient representation of \emph{context-specific independent} \cite{bout96} functions (CSI) where node $X$ is independent of nodes $W$ and $V$ given the context $u$ where node ($U=true$). Arithmetic operations can be performed on these functions returning a function value at the leaves; examples include addition
($\oplus$), subtraction ($\ominus$), multiplication ($\otimes$),
division ($\oslash$), $\min(\cdot,\cdot)$ and $\max(\cdot,\cdot)$ \cite{bahar93add}. 

Parameterized ADDs (PADDs) are an extension of ADDs that allow
for a compact representation of functions from $\{0,1\}^{n} \rightarrow
\mathbb{E}$, where $\mathbb{E}$ is the space of expressions
parameterized by $\vec{p}$. Formal definitions of our XADD are similar to that of PADD \cite{spuddip}.

Extended ADDs (XADDs) allow representing continuous variables in a decision diagram in the function representation of $\mathbb{R}^{n+m} \rightarrow \mathbb{R}$ over case statements. Each leaf in an XADD represents a multi-variate arbitrary function from the real-value domain and each decision node can be an equality, dis-equality or inequality on the multi-variate domain which is more expressive than the ADD boolean decisions. The branches are true/false depending on the value of each decision node. This compact representation will not require truth tables like ADDs or BDDs as it is more expressive in allowing infinitely many real values for each decision.

We next formally define the XADD operations and algorithms required to support all case operations of SDP as well as pruning algorithms to make this representation even more efficient. 

\subsection{Formal Definition and Operations}

An XADD allows polynomials at the leaves and decisions instead of a single real-value.  According to the set of continuous variables in an XADD, $\theta = \lbrace x_1,x_2,...,x_n\rbrace$ and the set of constants $c_i (0 \leq i \leq p)$ each leaf can be canonically defined as:
\begin{equation}
c_0+\sum_{i} c_i \prod_j \theta{ij}
\nonumber
\end{equation}
where $ 1 \leq j \leq n $. Each decision node is a polynomial inequality. Formally an XADD is defined using a BNF grammer: 
\begin{equation}
\nonumber
\begin{array}{lll}
F&::=&  \mathit{Poly}  \hspace{1mm} \vert  \hspace{1mm} \mathrm{if} (F^{var})\ \mathrm{then}\ F_h\ \mathrm{else}\ F_l \\
F^{var}&::=& ( \mathit{Poly} \leq 0)  \hspace{1mm}\vert \hspace{1mm}
 ( \mathit{Poly} \geq 0) \vert \hspace{1mm} B\\
\mathit{Poly}&::=&c_0+\sum_{i} c_i \prod_j \theta_{ij}\\
B&::=& 0|1
\end{array}
\end{equation}
An XADD node $F$ can either be a leaf $\mathit{Poly}$ with a polynomial value or a decision node $F^{var}$ with two branches $F_h$ and $F_l$ which are both of the non-terminal type $F$. The decision node $F^{var}$ associated with a single variable $var$ can be a polynomial inequality or a boolean decision $B = \lbrace b_1,b_2,...,b_m \rbrace$ where each boolean variable $b_k \in \lbrace 0,1 \rbrace$. If $F_h$ is taken the value of the decision node $F^{var}$ is true and if $F_l$ is taken the negation of the decision node $\neg F^{var}$ is set to true.\footnote {Note we assume continuous functions; if a function has the same values on a boundary point (equality), we allow only one of the $\leq, \geq$ at the boundary point. This continuous property allows us to replace $<$ ( $>$ ) with $\leq$($\geq$).} 

The value returned by a function $f$ represented as an XADD ($F$) containing (a subset of) discrete and continuous variables $\{ b_1,\cdots,b_m, x_1,\cdots, x_n \}$ with variable assignments $\rho \in \lbrace \lbrace\mathit{true}, \mathit{false} \rbrace^m,\mathbb{R}^n \rbrace$ can be defined recursively by:
\begin{equation*}
Val(F,\rho) = \left\{
\begin{array}{lll}
\textrm{if } F=\mathit{Poly} :&\mathit{Poly}\\
\textrm{if } F = \mathit{Poly} \rho(F^{var})=\mathit{true}: & \mathit{Val} (F_h,\rho)\\
\textrm{if } F = \mathit{Poly}  \rho(F^{var})=\mathit{false}: & \mathit{Val} (F_l,\rho)\\
\end{array} \right. 
\end{equation*}

This recursive definition of $Val(F,\rho)$ reflects the structural
evaluation of $F$ by starting at its root node and following
the branch at each decision node corresponding to the decisions taken
in $F^{var}$ --- continuing until a leaf node is reached,
which is then returned as $Val(F,\rho)$. The diagram on the right of Figure~\ref{fig:bdd_add_xadd} demonstrates the polynomial leaves and the decision node inequalities which branch to true/false depending on the decision value.

As with any tree-like structure, unreachable branches are bound to appear in XADDs, thus next we define pruning algorithms for more efficient results using the following definitions:

\begin{mydef}(\textbf{Function representation}):
A multi-variate function of booleans and real values $ f:\mathbb{B}^m \times \mathbb{R}^n \rightarrow \mathbb{R}$ denoted by $f(\vec{b},\vec{x})$  represents an XADD ($f_{XADD}$) defined on the class of piecewise formulas (case statements).
\end{mydef}

\begin{mydef}(\textbf{Path}):
A path $p$ in $f_{XADD}$ is a sequence of the pair $\rho=(F^{var},\mathit{dec})$ where each node  $F^{var}$ has a unique id and each decision assignment $\mathit{dec} \in \lbrace \mathit{true},\mathit{false}\rbrace$ represents the ( true or false) branch that node $F^{var}$ has followed. Note that the root node $\rho_0=(0,null)$ has a null decision assignment. A path is generally defined as a finite subset (in a sequence) of all possible pairs in $f_{XADD}$.
\begin{equation*}
p_j \subset \lbrace (F_1^{var},\mathit{dec_1}), (F_2^{var},\mathit{dec_2}), \cdots, (F_{n+m}^{var},\mathit{dec_{n+m}}) \rbrace
\end{equation*}
\end{mydef}
where $1 \leq j \leq (n+m)$ is the path number and the last pair on a given path $p_k$ is defined as the \textit{end-node} : $\eta(p_k)$ of that path.

\begin{mydef}(\textbf{Formula}):
The set of all paths for a node $F^{var}$ is defined as all paths $\lbrace p_1 \cdots p_k \rbrace$ ($1 \leq k \leq (n+m)$) such that the end-node of these paths are equal to node $F^{var}$ that is: $\eta(p_k) = F^{var}$. A formula on this node $\psi_{var} $ is defined as this finite set of paths $\psi_{var} = \lbrace p_k | \eta(p_k) = F^{var}, 1 \leq k \leq (n+m) \rbrace$. 
\end{mydef}

\vspace{10mm}
\begin{figure}[h]
\centering

\includegraphics[width=0.8\linewidth]{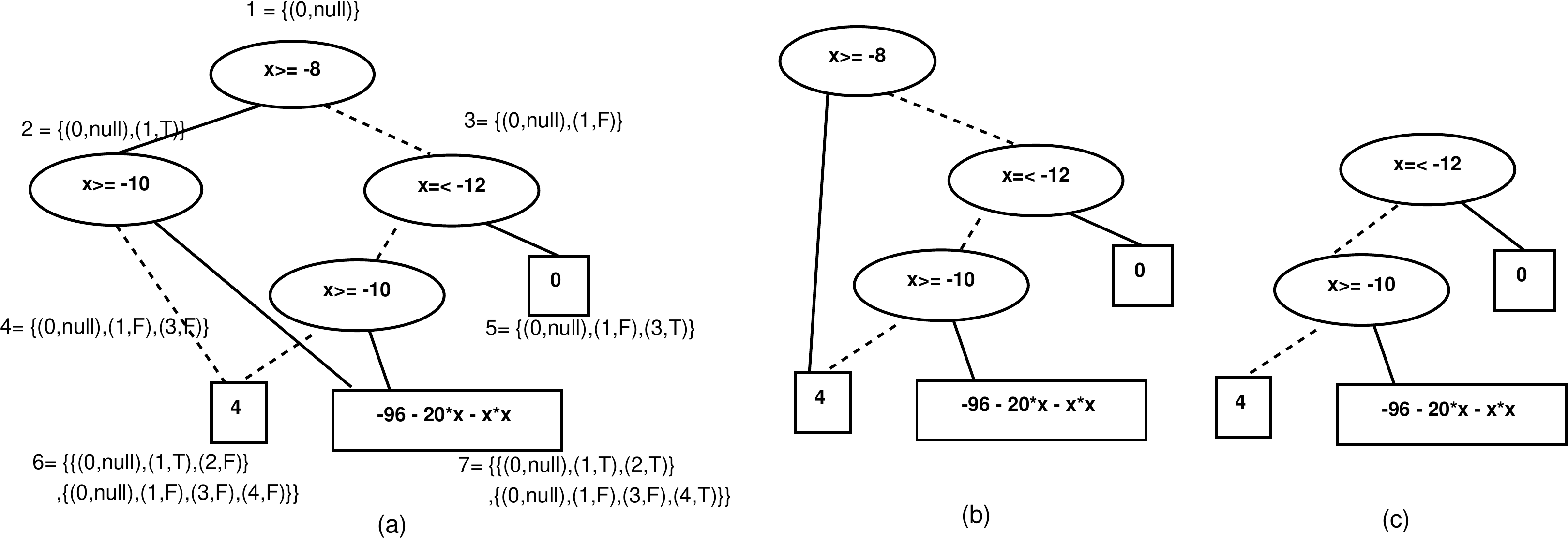}

\caption{\footnotesize (a) Path and formula definitions on each node in the XADD. A path is defined as a sequence of the tuple $(F^{var},\mathit{dec})$ and a formula is a set of paths. (b) Pruning inconsistent nodes of (a). (c) Pruning redundant nodes of (b)}
\label{fig:path_formula}

\end{figure}

To define any node $F^{var}$ in $f_{XADD}$ logically we have the following:
\begin{align}
\psi_{var} &= \bigvee_{p_j \in p_{var}} \left( \bigwedge_{\rho_i \in p_j} \rho_i \mid \eta(p_j) = F^{var} \right) \nonumber \\
& \hspace{5mm} , \nonumber \\
\rho_i &= \begin{cases}
  \text{dec}_i = \mathit{true}: & F_i^{var} \\ 
  \text{dec}_i = \mathit{false}: & \neg F_i^{var} \\ 
\end{cases} \label{eq:rho_i}
\end{align}
Figure~\ref{fig:path_formula} (a) shows a simple XADD with all paths and formulas determined on each node.
\begin{mydef}(\textbf{Node ordering}):
A node $F_i^{var}$ in $f_{XADD}$ is defined before node $F_j^{var}$ if it appears before this node in a path $p_s$ containing both nodes $\rho_i \in p_s, \rho_j \in p_s$ and has an ordering such that $i<j$. $F_j^{var}$ can also be named the child of node $F_i^{var}$. A parent node $F_k^{var}$ is defined for node  $F_i^{var}$ if node  $F_i^{var}$ appears after node $F_k^{var}$ in a path $p_s$ containing both nodes $\rho_i \in p_s, \rho_k \in p_s$ and has an ordering such that $k<i$. In Figure~\ref{fig:path_formula} (a) node $1$ is the parent of node $3$ and nodes $4,5$ are the children of node $C$.
\end{mydef}

\begin{mydef}(\textbf{Inconsistent node}):
A node $F_i^{var}$ in $f_{XADD}$ is inconsistent if it violates any of the constraints in its parent decision node $F_k^{var}$ in any of the path defined in the formula $\psi_{i}$ over this node. In mathematical terms this is equal to the following: 
 \begin{align}
\exists p_j \in \psi_i, (\rho_i \in p_j): \quad & \exists \rho_k \in p_j, k < i, \nonumber \\
& \phi_i \Rightarrow \phi_k = \textit{true} \nonumber \\
& \longleftrightarrow F_i^{var} = \text{inconsistent}. \label{inconsist}
\end{align}
where  $\phi_i$ is the logical constraints of node $F_i^{var}$ as defined in case statements. Figure~\ref{fig:path_formula} (b) prunes (a) of the inconsistent node 2.
\end{mydef}

\begin{mydef}(\textbf{Redundant node}):
A hub $F_i^{var}$ in $f_{XADD}$ is excess in the event that its requirements can be tended to involving any of the imperatives in its kid choice hub $F_j^{var}$ in any of the way characterized in the equation $\psi_{j}$ over this hub. In numerical terms this is equivalent to the accompanying:
\begin{align}
\exists p_k \in \psi_j, (\rho_j \in p_k): \quad & \exists \rho_i \in p_k, \, i < j, \nonumber \\
& (\phi_j \Rightarrow \phi_i) \vee (\neg \phi_j \Rightarrow \phi_i) = \textit{true} \nonumber \\
& \longleftrightarrow F_i^{var} = \text{redundant}. \label{redundant}
\end{align}

Figure~\ref{fig:path_formula} (c) prunes (b) of the repetitive hub 1.
\end{mydef}

For any capability from $\mathbb{R}^{n+m} \rightarrow \mathbb{R}$, we next
depict how a decreased XADD can be built from an inconsistent arranged choice graph.
All calculations that we will deﬁne in the accompanying segments depend on the partner capability \emph{getNode} in Calculation \ref{algGetNode} , which returns a more conservative portrayal of a solitary interior choice hub.

The calculation \emph{ReduceXADD} permits the development of a conservative XADD portrayal from an erratic arranged choice chart with polynomial leaves and polynomial imbalances as the choice hubs. Algorithm~\ref{algReduceXADD} is characterized by the accompanying definition:
\incmargin{1em}
\linesnumbered
\begin{algorithm}[h]
\SetKwFunction{getNode}{{\sc GetNode}}
\SetKwFunction{reduce}{{\sc ReduceXADD}}
\SetKwInOut{Input}{input}
\SetKwInOut{Output}{output}

\Input{$F$ (root node id for an arbitrary ordered decision diagram)}
\Output{$F_r$ (root node id for reduced XADD)}
\BlankLine
\Begin{
   //if terminal node, return canonical terminal node\\
   \If{F is terminal node}
   {
   \Return{canonical terminal node for polynomial of $F$}\;
   }
   //use recursion to reduce sub diagrams\\
   \If{$F \rightarrow F_r$ is not in ReduceCache}
   {
    $F_h$ = \reduce{$F_h$}\;
    $F_l$ = \reduce{$F_l$}\;
    //get a canonical internal node id\\
    $F_r$ = \getNode{$F^\mathit{var}$, $F_h$, $F_l$}\;
    insert $F \rightarrow F_r$ in ReduceCache\;
   } 
   \Return{$F_r$}\;
}
\caption{{\sc ReduceXADD}(F)  \label{algReduceXADD}}
\end{algorithm}
\decmargin{1em}
\begin{algorithm}[h]
\SetKwInOut{Input}{input}
\SetKwInOut{Output}{output}

\Input{$\langle \mathit{var}, F_h, F_l \rangle $ (variable and true and false branches node ids for internal node)}
\Output{$F_r$ (canonical internal node id)}
\BlankLine
\Begin{
   //redundant branches\\
   \If{$F_l = F_h$}
   {
      \Return{$F_l$}\;
   }
   //check if the node exists previously\\
   \If{$\langle \mathit{var}, F_h, F_l \rangle \rightarrow \mathit{id}$ is not
   in NodeCache}
   {id\ = new unallocated id\;  
    insert $\langle \mathit{var}, F_h, F_l \rangle \rightarrow \mathit{id}$ in NodeCache\;
   }
   \Return{id}\;
}
\caption{{\sc GetNode}($\langle \mathit{var}, F_h, F_l
\rangle $)  \label{algGetNode}}
\end{algorithm}

\textbf{Definition:} A capability chart \( G \) is decreased on the off chance that it contains no vertex \( v \) with \( \text{low}(v) = \text{high}(v) \), nor does it contain particular vertices \( v \) and \( v' \) to such an extent that the subgraphs established at \( v \) and \( v' \) are isomorphic.

This calculation recursively builds a decreased XADD from the base up. Interior hubs are addressed as $\langle F^{\mathit{var}}, F_h, F_l \rangle$, where $F^{\mathit{var}}$ is the variable name, and $F_h$ and $F_l$ are the valid and misleading branch hub ids, individually. Diminished hubs are put away in the \emph{ReduceCache} table.
Utilizing the capability \emph{GetNode} (Algorithm~\ref{algGetNode}) any excess choice tests are eliminated. This capability stores an extraordinary id for every hub in the \emph{NodeCache} table.

\emph{ReduceCache} guarantees that every hub is visited once and a special decreased hub is produced in the last chart. Subsequently \emph{ReduceXADD} has straight running reality as indicated by the size of the info diagram.
\vspace{10mm}
\begin{figure}[h]
\centering

\includegraphics[width=0.65\linewidth]{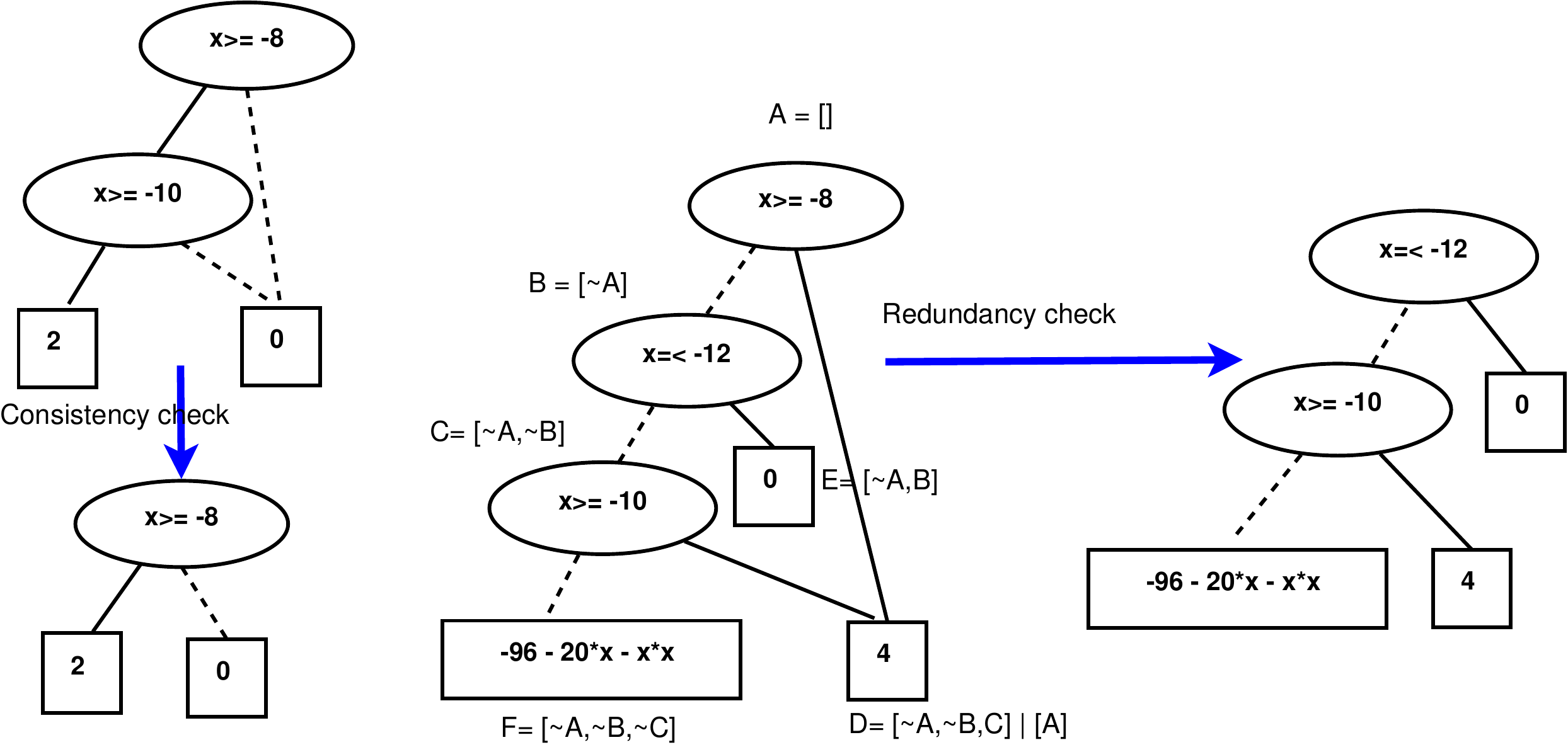}

\hspace{10mm}
\includegraphics[width=0.2\linewidth]{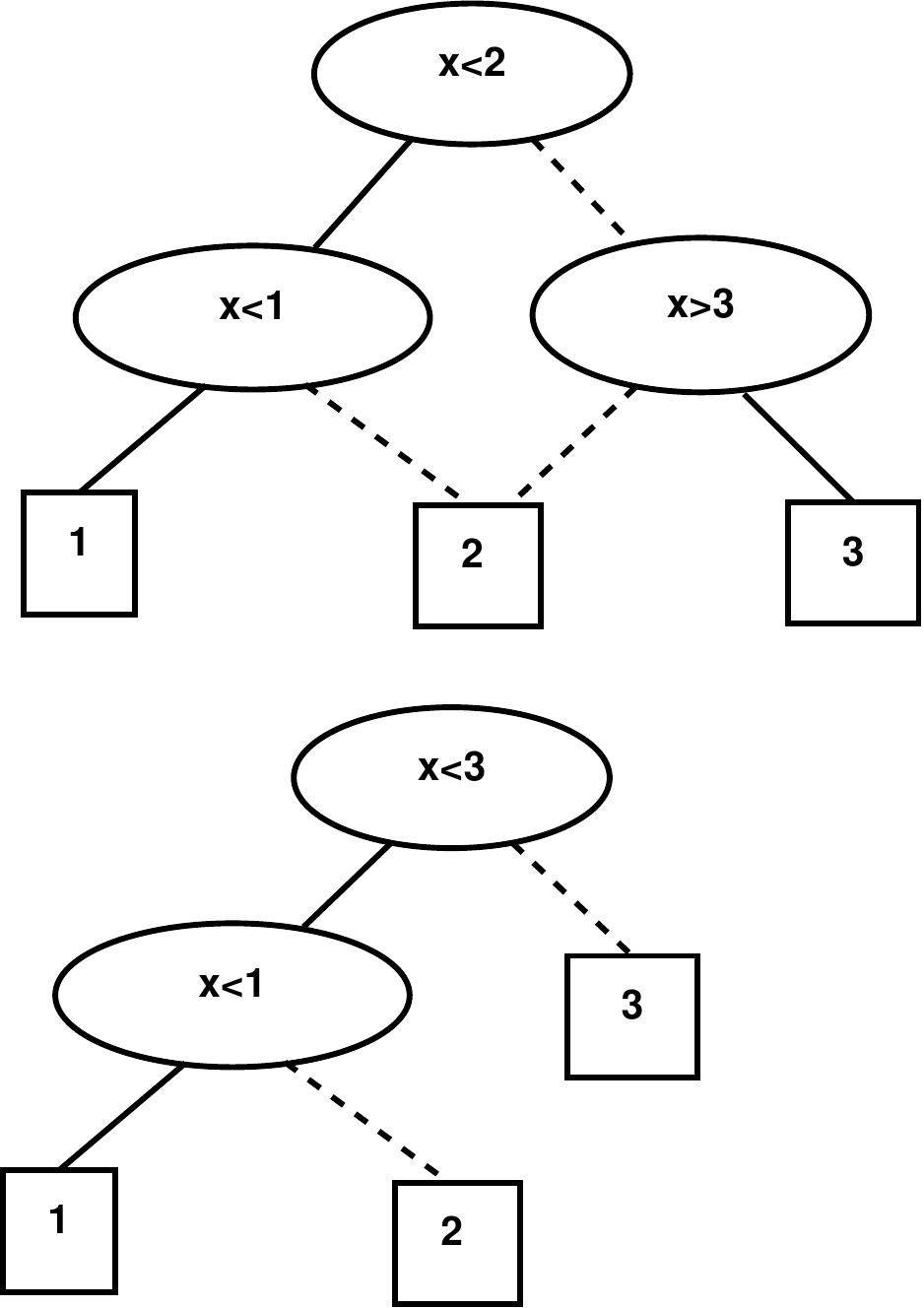}
\vspace{-2mm}

\caption{\footnotesize (Left) Using pruning algorithms for inconsistency and redundancy. Top-left figure is reduced to the bottom-left figure using the LP-Solver which recognizes the parent-child relation. The middle figure is reduced to the right-most figure using the SAT-Solver. The paths of each node in the middle graph are demonstrated. These paths along with the child-parent implications in the KB can reduce the redundant node A. (Right) A counterexample for an XADD that is not canonical even after applying consistency and redundancy checking. The bottom diagram is the true value of the XADD which can not be derived after pruning.}
\label{fig:canonical}

\end{figure}

We next present two successive algorithms, Algorithm~\ref{algPrune} for removing inconsistent nodes and Algorithm~\ref{algRedundant} for removing redundant nodes. 
\incmargin{1em}
\linesnumbered
\begin{algorithm}[h]
\SetKwFunction{getNode}{{\sc GetNode}}
\SetKwFunction{testImplied}{{\sc TestImplied}}
\SetKwFunction{prune}{{\sc PruneInconsistent}}
\SetKwInOut{Input}{input}
\SetKwInOut{Output}{output}

\Input{$F$ (root node id for an inconsistent XADD, Decision label, Decision value)}
\Output{$F_r$ (root node id for a consistent XADD)}
\BlankLine
\Begin{
   //if terminal node, return the value\\
   \If{F is terminal node}
   {
    \Return{canonical terminal node for polynomial of $F$}\;
   }
   //else if internal node, find all implications consider any parent in XADD\\
   
    	//if (parent$\Rightarrow$ child) remove child\\
    	 \If{\testImplied{$\mathit{decLabel}$,$\mathit{decision}$,$F^\mathit{var}$}}
    			{$F_r$ = \prune{$F_h$,$\mathit{decLabel}$,$\mathit{decision}$}\;}
    	//if (parent$\Rightarrow$ $\neg$child) remove child\\
    	 \If{!\testImplied{$\mathit{decLabel}$,$\mathit{decision}$,$- F^\mathit{var}$}}
    			{$F_r$ = \prune{$F_l$,$\mathit{decLabel}$,$\mathit{decision}$}\;}
  //if result of LP-solver is null\\
   \Else
   {	
     insert $F^{var} \rightarrow \mathit{decLabel}$\;
     insert $\mathit{false} \rightarrow \mathit{decision}$\;
    $F_r^l=$  \prune{$F_l$,$\mathit{decLabel}$,$\mathit{decision}$}\;
     insert $\mathit{true} \rightarrow \mathit{decision}$\;
    $F_r^=$  \prune{$F_h$,$\mathit{decLabel}$,$\mathit{decision}$}\;	 
     insert $\mathit{null} \rightarrow \mathit{decision}$\;
     delete $F^{var} \rightarrow \mathit{decLabel}$\;

    $F_r=$ \getNode{$var,F_r^h,F_r^l$}\;
	}
   \Return{$F_r$}\;
}
\caption{{\sc PruneInconsistent}(F,$\mathit{decLabel}, \mathit{decision}$)  \label{algPrune}}
\end{algorithm}
\decmargin{1em}

Given any XADD with potential inconsistent nodes, the output of \emph{algPrune} (Algorithm~\ref{algPrune}) is a reduced XADD with canonical leaves, linear decisions and no inconsistent nodes. 
In lines 7 and 11 we test the current decision node with all previous decisions for implications in \emph{TestImplied}. 
This process is performed using the definition of constraints for the LP-solver. Once a node is returned inconsistent (returns true or false from \emph{TestImplied}) the high or low branch is returned without referring to the current node in the final XADD. 

While any (parent(A) $\Rightarrow$ child(B)) relation is returned for consistency checking, at the same time this function stores other implications in the a related cache ($\mathit{hmImplications}$). Since $A \Rightarrow B$ is already tested, we cache the following child to parent implications: (i) $B \Rightarrow A$ (2)$\neg B \Rightarrow A$ ($\neg B \Rightarrow \neg A$ is already defined by $A \Rightarrow B$).We add these to the knowledge-base (KB) of the SAT-Solver before performing the next recursive algorithm. 

\incmargin{1em}
\linesnumbered
\begin{algorithm}[h]
\SetKwFunction{getNode}{{\sc GetNode}}
\SetKwFunction{sat}{{\sc Sat-Test}}
\SetKwFunction{prune}{{\sc PruneRedundancy}}
\SetKwInOut{Input}{input}
\SetKwInOut{Output}{output}
\Input{$X$ (root node id for a consistent XADD), $KB$ (child-parent implications), $F$ (formulas for each node)}
\Output{$X_r$ (root node id for a consistent and redundant XADD)}
\BlankLine
\Begin{
	//current node is X, keep its path\\
	$\mathit{Path}:=F(X)$\;
   \If{F is terminal node}
   {
    \Return{canonical terminal node for polynomial of $F$}\;
   }
   //else if internal node, add paths to low and high branch\\
   \ForEach{$\mathit{path} \in$  $\mathit{Path}$ }
   {
    	//add high and low branch of internal node to current path\\
    	insert $X_r \rightarrow  \mathit{path};$
    	add $\mathit{path} \rightarrow F$	\; 
    	insert $X_l \rightarrow  \mathit{path};$
    	add $\mathit{path} \rightarrow F$	\;
    	
    }	
    $X_l=$  \prune{$X_l,KB,F$}\;
    $X_h^=$  \prune{$X_h,KB,F$}\;	 
   
   //reached the lowest decision node, perform satisfiability test\\
    \If{\sat{$X$,$F(X_l),KB,T$}}
    	{ \Return{$X_h$}\;}
    \If{\sat{$X$,$F(X_h),KB,F$}}
    	{ \Return{$X_l$}\;}
   
    $X_r=$ \getNode{$var,X_h,X_l$}\;
	
   \Return{$X_r$}\;
}
\caption{{\sc PruneRedundancy}(X,KB,F)  \label{algRedundant}}
\end{algorithm}
\decmargin{1em}
Lines 14--15 define the decision labels and their values ($\mathit{true,false}$) for the current node and add them to the set of labels and decisions which are input to the algorithm. The final result is returned from the \emph{GetNode} algorithm. As an example consider the two cases on the left side of the Figure~\ref{fig:canonical}. In the upper diagram,  $x \geq -8  \Rightarrow x \geq -10 $ is reduced to $x\geq-8$ since it covers all the state space in the lower diagram.   

For redundancy pruning an equivalence testing approach  is proposed using a SAT-solver. Consider the middle example in Figure~\ref{fig:canonical} where the child node(C) implies the parent node(A). We define the recursive algorithm required to prune all child-parent implications from an XADD using this example. We introduce Algorithm \ref{algRedundant} to remove redundant nodes. 
The input to this algorithm is an XADD where each node $F^{var}$ is marked with its formula $\psi_{var}$ (set of paths) as in the middle diagram of Figure~\ref{fig:canonical}. 

The recursive property of the algorithm traverses backwards from the leaf nodes up to the root node. Every time a decision node is reached, the true and false branches of that node are put to a satisfiability test (\emph{Sat-Test}) to see if that branch can be eliminated from the tree due to some child-parent implication.
We test the formulas $\psi_{A=F/T} \Longleftrightarrow \psi$ for Figure~\ref{fig:canonical} this is equal to the following satisfiability test on the false branch of node A: 
\begin{equation*}
(B \Rightarrow A) \models ((A \vee (\neg A \wedge \neg B \wedge C))\Longleftrightarrow(\mathit{false} \vee (\mathit{true} \wedge \neg B \wedge C)))
\end{equation*}

We use a SAT-solver (minisat ~\cite{minisat}) to entail this sentence where the KB here contains the child-parent implication $C \Rightarrow \neg A$  and since the result is true, we can prune the true branch returning node B represented in the right-hand diagram of  Figure~\ref{fig:canonical}. A similar approach is performed for the both branches in lines 14--17. The returned value is the true or false branch if the parent node can be pruned else a new node is built using line 19 of Algorithm \ref{algGetNode}.

For BDDs and ADDs proof of canonicity can be defined according to \cite{bryant}. In these two cases for any function $f(x_1,\cdots, x_n)$ and a fixed variable ordering over $x_1,\cdots, x_n$, a reduced BDD or ADD is defined as the minimally sized ordered decision diagram representation of a function $f$. This proves that there is a unique canonical BDD or ADD representation for every function from $\{0,1\}^{n} \rightarrow \{0,1\}$ for BDDs and $\{0,1\}^{n} \rightarrow \mathbb{R}$ for ADDs. 

Unlike BDDs and ADDs proving canonicity for XADDs is not applicable. As a counterexample consider the simple XADD in the top right diagram of Figure~\ref{fig:canonical} where the values as below: 
\begin{align*}
	\begin{cases}
		x \leq 1 :& 1  \\
		x \geq 3 :&3\\
		1<x<3 :& 2\\
	\end{cases}
\end{align*}
As the values suggest, there is no need to branch on $x<2$ in this function thus it can be removed from the tree. According to the pruning techniques mentioned, this XADD is neither inconsistent nor redundant. To remove the node, a reordering of the decision nodes is required which will effect the structure of the tree, changing it to a new XADD. For this reason, although consistency and redundancy checking takes care of most unwanted branches, there is no guarantee that a given XADD is minimal after applying the two pruning techniques. The bottom right diagram of Figure~\ref{fig:canonical} is the true value this XADD holds which is a different XADD.

In principle exact SDP solutions can be obtained for arbitrary symbolic functions, we restrict XADDs to use polynomial functions only.
Having proved the main properties of an XADD, next we present the operations and algorithms required for SDP for XADDs.

\subsection{XADD algorithms and operations}

In this section we review the symbolic operations required to perform SVI using the XADD structure. This is mainly categorized into unary and binary operations.   

\begin{algorithm}[h]
\SetKwFunction{getCanonicalNode}{{\sc GetCanonicalNode}}
\SetKwFunction{reduce}{{\sc Reorder}}
\SetKwInOut{Input}{input}
\SetKwInOut{Output}{output}

\Input{$F$ (root node for possibly unordered XADD)}
\Output{$F_r$ (root node for an ordered XADD)}
\BlankLine
\Begin{
   //if terminal node, return canonical terminal node\\
   \If{F is terminal node}
   {
   \Return{canonical terminal node for polynomial of $F$}\;
   }
   //else nodes have a $\mathit{true}$ \& $\mathit{false}$ branch and $\mathit{var}$ id\\
   \If{$F \rightarrow F_r$ is not in Cache}
   {
    $F_{\mathit{true}}$ = \reduce{$F_{\mathit{true}}$} $\otimes \; \mathbb{I}[F_\mathit{var}]$ \;
    $F_{\mathit{false}}$ = \reduce{$F_{\mathit{false}}$} $\otimes \; \mathbb{I}[\neg F_\mathit{var}]$\;
    $F_r = F_{\mathit{true}} \oplus F_{\mathit{false}}$\;
    insert $F \rightarrow F_r$ in Cache\;
   } 
   \Return{$F_r$}\;
}
\caption{{\sc Reorder}(F)  \label{alg:reorder}}
\end{algorithm}

\subsubsection{Unary XADD Operations}
According to the previous section on SDP unary operations, scalar multipication $c.f$ and negation $- f$ on a function results in a function which can simply be represented as an XADD. Apart from this restriction, substitution and marginalization of the $\delta$ function are also unary operations that can be applied to XADDs are explained below.  

Restriction of a variable $x_i$ in an XADD ($F$) to some formula $\phi$ is performed by appending $\phi$ to each of the decision nodes while leaves are not affected.
For a binary variable restriction to a single variable $x_i$ is equal to taking the true or false branch according to that variable 
($F|_{x_i=true}$) or ($F|_{x_i=false}$). 
This operation can also be used to $\texttt{marginalize}$ or $\texttt{sum out}$ boolean variables in a CA-HMDP. 
$\sum_{x_i \in X_i}$ eliminates a variable $x_i$ from an XADD and is computed as the sum of the \emph{true} and \emph{false} restricted functions ($F|_{x_i=true} \oplus \ F|_{x_i=false}$). We omit the restriction and marginalization algorithms for XADDs since they are identical to the same operations for ADDs. 

Substitution for a given function $f$ is performed by applying $\sigma$ the set of variables and their substitutions to each inequality operand such that $\phi_i\sigma: f_i\sigma$. The substitution operand effects both leaves and decision nodes and changes them according to the variable substitute. 

Decisions become unordered when substituted and also when we perform maximum or minimum explained in the next section. A reorder algorithm has to be applied to the result of the substitution operand. As Algorithm~\ref{alg:reorder} shows, we recursively apply the binary ADD operations of $\otimes$ and $\oplus$ to decision nodes to reorder the XADD after a substitution. 

As for the integration of the $\delta$-function on variable $x$ we require computing $\int_{x} \delta [ x - g(\vec{x})]fdx$. This triggers the substitution $f \lbrace x/ g(\vec{x})\rbrace$ on $f$ as defined above.  

For the continuous maximization $max_y$ each XADD path from root to leaf node is considered a single case partition with conjunctive constraints, and maximization is performed at each leaf subject to these constraints and all path maximums are then accumulated using the \textit{casemax} operation for the final result. Note that in general continuous maximization is a multi-variate operation but according to the previous section, it can be decomposed into multiple univariate operations. 

\vspace{10mm}
\begin{figure}[h]

\includegraphics[width=0.97 \linewidth]{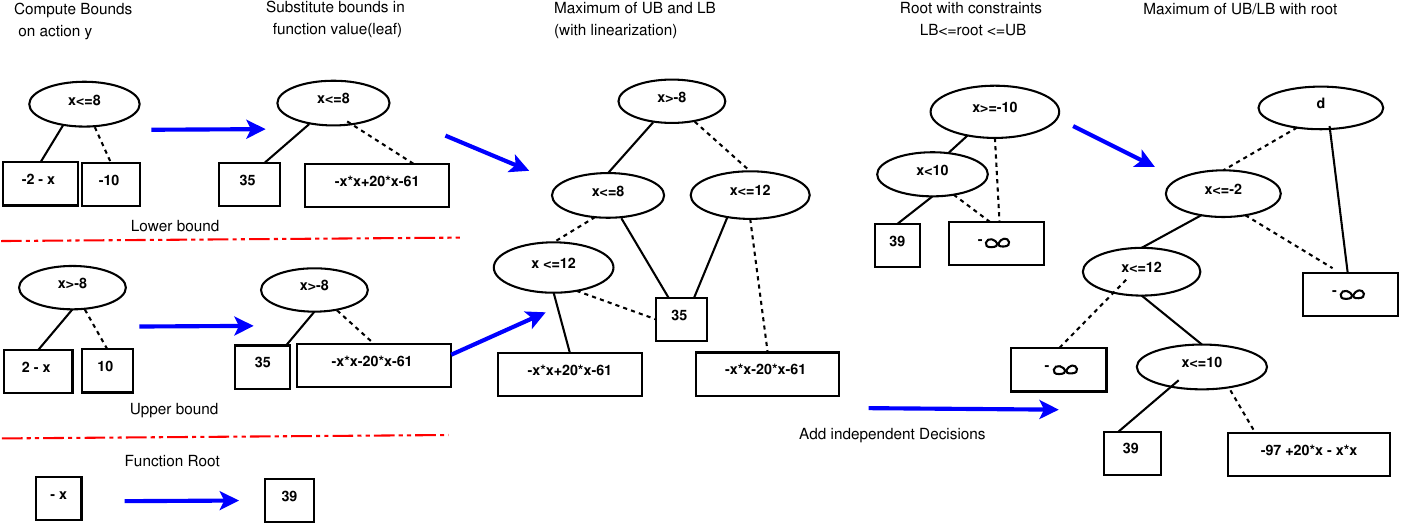}

\caption{\footnotesize Representing one step of Continuous Maximization algorithm using XADDs. This partition has conditions: 
$\phi_i ( x,d,y ) \equiv \neg d \wedge ( x>2 ) \wedge ( -10 < y < 10) \wedge ( -2 < x+y < 2 )$  
and value $f_i ( x,y ) = 39- y^2 -x^2 - 2a \times x$. The upper and lower bounds of action $y$ and the root is represented by XADDs and substituted inside the leaf node. The maximum of the upper and lower bound is then represented after linearizing the result. The final result takes the maximum of this and the root considering the constraints and independent decisions.}
\label{fig:xadd_max}
\vspace{-8mm}
\end{figure}

A univariate continuous maximization algorithm has been presented in the previous section. 
Starting at the root node, if the current node is a decision node, the algorithm is recursively calls on the low and high branches of this node. If the current node is a leaf node, it processes the leaf according to \emph{Continuous Maximization}. The input to this algorithm is the leaf node and the decisions leading to this leaf node i.e., partitioning of the state-action space. The algorithm finds the maximum value according to XADDs used for the upper and lower bounds as well as the function roots. The maximum XADD resulting from these three XADDs along with applying the independent constraints define the final result. A single step of this algorithm is presented in Figure~\ref{fig:xadd_max}. The input to this step is the leaf node obtained after the regression step of Algorithm \ref{alg:regress}. Also as an input are the decisions leading to this leaf (high or low branches).  The XADD representation is used to represent all intermediate results as well as the final result. 

In order to solve the problem of non-linear domains we require a continuous maximization of the continuous action variable $y$ over a non-linear function. This maximization is performed symbolically by incooperating the non-linear terms at the nodes inside decision nodes. This does not effect our symbolic solution but in order to prune the resulting value function using a linear programming solver, we require linear decisions. 

The linearize algorithm (Algorithm \ref{algLinearize}) presented here is of a recursive nature similar to other XADD algorithms. For each node in the XADD starting at the root node, the algorithm linearizes the decision node using \emph{GetRoots} to find the roots of the non-linear function. It then returns the decision nodes replacing that single non-linear node and  iterates on the low and high branches of the decision node until it returns all the leaves 			
\incmargin{1em}
\linesnumbered
\begin{algorithm}[h]
\SetKwFunction{getRoots}{{\sc GetRoots}}
\SetKwFunction{reduce}{{\sc ReduceLinearize}}
\SetKwInOut{Input}{input}
\SetKwInOut{Output}{output}

\Input{$F$ (root node id for an arbitrary ordered decision diagram)}
\Output{$F_r$ (root node id for linearized XADD)}
\BlankLine
\Begin{
   //if terminal node, return canonical terminal node\\
   \If{F is terminal node}
   {
   \Return{canonical terminal node for polynomial of $F$}\;
   }
    \If{$F \rightarrow F_r$ is not in ReduceCache}
   {
   //use recursion to reduce sub diagrams\\
    $F_h$ = \reduce{$F_h$}\;
    $F_l$ = \reduce{$F_l$}\;
    //get a linearized internal node id\\
    $F_r$ = \getRoots{$F^\mathit{var}$, $F_h$, $F_l$}\;
    insert $F \rightarrow F_r$ in ReduceCache\;
   } 
   \Return{$F_r$}\;
}
\caption{{\sc ReduceLinearize}(F)  \label{algLinearize}}
\end{algorithm}
\decmargin{1em}

\subsubsection{Binary XADD Operations}

For \emph{all} binary operations, the function \emph{Apply}($F_1,F_2,\mathit{op}$) computes the resulting XADD. 

Two canonical XADD operands $F_1$ and $F_2$ and a binary operator $\mathit{op} \in \{ \oplus, \ominus, \otimes , \max , \min \} $ are the input to the \emph{Apply} algorithm. The output result is a canonical XADD $F_r$. 
If the result of \emph{Apply}($F_1,F_2,\mathit{op}$) is a quick computation as in line 3, 
it can be immediately returned. Else it checks the \emph{ApplyCache} in line 6 for any previously stored apply result.  If there is not a cache hit, the earliest variable in the ordering to branch is chosen according to \emph{ChooseEarliestVar}
Two recursive \emph{Apply} calls are then made on the branches of this variable to compute $F_l$ and $F_h$. \emph{GetNode} checks for any redundancy in line 23 before storing it in the cache and returning the resulting XADD. We cover these steps in-depth in the following sections.

\begin{algorithm}[h]
\SetKwInOut{Input}{input}
\SetKwInOut{Output}{output}

\Input{$F_1$ (root node id for operand 1),\\
$F_2$ (root node id for operand 2)}
\Output{$var$ (selected variable to branch)}
\Begin{
   //select the variable to branch based on the order criterion\\
   \eIf{$F_1$ is a non-terminal node}
    {
        \eIf{$F_2$ is a non-terminal node}
         {
	    \eIf{$F_1^{var}$ comes before $F_2^{var}$}
	    {
                $\mathit{var}=F_1^{var}$\;
	    }
	    {
	        $\mathit{var}=F_2^{var}$\;
	    }
         }
         {
           $\mathit{var}=F_1^{var}$\;
         }
    }
    {
     $\mathit{var}=F_2^{var}$\;
    }
   \Return{$\mathit{var}$}\;
}
\caption{{\sc ChooseEarliestVar}($F_1,F_2$)  \label{algChooseVarBranch}}
\end{algorithm}

\begin{table}[h]
\begin{center}
    \begin{tabular}{|l|l|l|}
	 \hline
	 Case number&Case operation & Return \\ \hline \hline
1&         $F_1\ \mathit{op}\ F_2; F_1=\mathit{Poly}_1; F_2=\mathit{Poly}_2$ & $\mathit{Poly}_1\ \mathit{op}\ \mathit{Poly}_2$ \\
	 \hline
2&          $F_1  \oplus F_2; F_2=0$ & $F_1$\\
	 \hline
3&          $F_1  \oplus F_2; F_1=0$ & $F_2$\\
	 \hline
4&         $F_1  \ominus F_2; F_2=0$ & $F_1$\\
	 \hline
5&         $F_1  \otimes F_2; F_2=1$ & $F_1$\\
	 \hline
6&          $F_1  \otimes F_2; F_1=1$ & $F_2$\\
	 \hline
7&         $F_1  \otimes F_2; F_2=0$ & 0\\
	 \hline
8&          $F_1  \otimes F_2; F_1=0$ &0\\
	 \hline
9&          $F_1  \oplus \infty $ & $\infty$\\
	 \hline
10&          $F_1  \otimes \infty ; F_1 \neq 0$ & $\infty$\\
	 \hline
11&          $\max (F_1  , +\infty)$ & $\infty$\\
	 \hline
12&          $\max (F_1, -\infty)$ & $F_1$\\
	 \hline
13&          $\max (F_1  , F_2)$ &pic\\
	 \hline
14&          $\min (F_1, F_2)$ &pic\\
	 \hline
13&	  other& $\mathit{null}$\\
         \hline
    \end{tabular}
  \caption{Input case and result for the method \emph{ComputeResult}
  for binary operations  $\oplus$,  $\ominus$ and $\otimes$ for XADDs.}
  \label{tab:ComputeResultXADD}
\end{center}
\vspace{-10mm}
\end{table}
\subsubsection{Apply Algorithm for binary operations of XADDs}


\textbf{Terminal computation}: 
The function \emph{ComputeResult}  determines if the result of a computation can be immediately computed without recursion. The entries denote a number of pruning optimizations that immediately return a node without recursion. For the discrete maximization (minimization) operation (entries 11--14) , for every two leaf nodes $f$ and $g$ an additional decision node $f > g$ ($f < g$) is introduced to represent the maximum(minimum). This may cause out-of-order decisions which can be solved by the reordering Algorithm~\ref{alg:reorder}. 

\textbf{Caching}:
If the result of 	\emph{ComputeResult} is empty, in the next step we check the  \emph{ApplyCache} for any previously computed operation using this set of operands and operations. To increase the chance of a match, all items stored in a cache are made canonical.

\textbf{Recursive computation}:
If a call to Apply is unable to immediately compute a result or reuse a previously cached computation, we must recursively compute the result. If both operands are constant terminal nodes the function \emph{ComputeResult} takes care of the result and for the rest of the cases one of the following conditions applies: 
\begin{itemize}
\item $F_1$ or $F_2$ is a constant terminal node or $F_1^{var} \neq F_2^{var} $: The high and low branch of the operand chosen by \emph{ChooseEarliestVar} (here $F_2$) is applied with the other operand $F_1$ and the result of these two \emph{Apply} functions are used with an $\mathit{if}$ decision node on $F_2^{var}$ for a canonical result:
\begin{align*}
F_h &= Apply (F_1 , F_{2,h},op) \\
F_l &= Apply (F_1 , F_{2,l},op)  \\
F_r &= GetNode (F_2^{var} , F_h,F_l) 
\end{align*}
\item $F_1$ and $F_2$ is constant nodes and $F_1^{var} = F_2^{var} $: Here branching is an $if$ statement on $F_1^{var} (=F_2^{var})$ with the true case applied to the high branches of $F_1$ and $F_2$ and the false case is defined as the result of the \emph{Apply} function on the low branches:
\begin{align*}
F_h &= Apply (F_{1,h} , F_{2,h},op)  \\
F_l &= Apply (F_{1,l} , F_{2,l},op)  \\
F_r &= GetNode (F_1^{var} , F_h,F_l)  \\
\end{align*}
\end{itemize}

Finally for the high and low branch of the final computation, two recursive calls are made to \emph{Apply} and the result is in a canonical form returned by \emph{GetNode}. 
Having defined the efficient representation of XADDs, next we show results from implementing the SVI algorithms using this structure.

\section{Experimental Results} \label{results}
\vspace*{-0.05in}

We implemented two versions of our proposed SVI algorithms using XADDs
--- one with the discrete setting and one with the continuous setting.

For CA-HMDPs we evaluated SVI on a didactic nonlinear
\MarsRover\ example and two problems from Operations Research (OR) \InventoryControl\ defined in the introduction and \WaterReservoir  all of which are described below. For comparison purposes the DA-HMDPs example domains are discretized by their action space. \footnote{All Java source code and a human/machine readable file format for all domains needed to reproduce
the results in this paper can be found online at
\texttt{http://code.google.com/p/xadd-inference}.}

\subsection{Domains}

\paragraph{\MarsRover}
A \MarsRover state comprises of its consistent position $x$ along a given course. In a given time step, the meanderer may move a consistent distance $y \in [-10,10]$. The wanderer accepts its most noteworthy compensation for snapping a photo at $x=0$, which quadratically diminishes to zero at the limits of the reach $x \in [-2,2]$. The meanderer will
naturally snap a photo when it begins a period step inside the reach $x \in [-2,2]$ and it just gets this prize once.

Utilizing boolean variable $b \in \{0,1\}$ to demonstrate assuming the image has
as of now been taken ($b=1$), $x'$ and $b'$ to indicate
post-activity state, and $R$ to indicate reward, we
express the \MarsRover\ CA-HMDP utilizing piecewise elements and prize:
\begin{align*} 
\hspace{-2.8mm} P(b'\sq=\sq1|x,b) & = 
\begin{cases}
b \lor (x \geq -2 \land x \leq 2): & \sqm 1.0\\
\neg b \land (x < -2 \lor x > 2):  & \sqm 0.0
\end{cases}  \\
\hspace{-2.8mm} P(x'|x,y) & = \delta \left( x' - \begin{cases}
y \geq -10 \land y \leq 10 : & \hspace{-2mm} x + y \\
y < -10 \lor y > 10 : & \hspace{-2mm} x
\end{cases}
\right)  \\
\hspace{-2.8mm} R(x,b) & = \begin{cases}
\neg b \land x \geq -2 \land x \leq 2 : & 4 - x^2 \\
b \lor x < -2 \lor x > 2 : & 0
\end{cases} 
\end{align*}
\begin{figure}[h]

\begin{minipage}[b]{0.8\linewidth}
\includegraphics[width=0.9\textwidth]{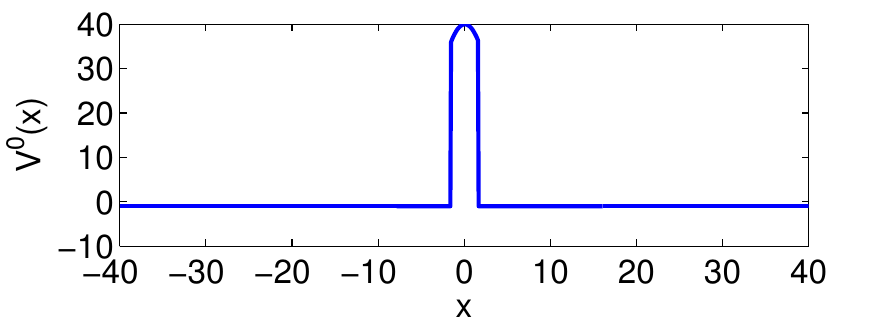}\\

\includegraphics[width=0.9\textwidth]{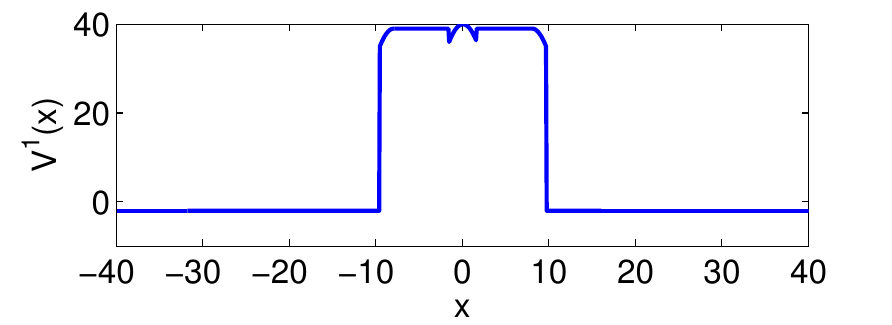}\\

\includegraphics[width=0.9\textwidth]{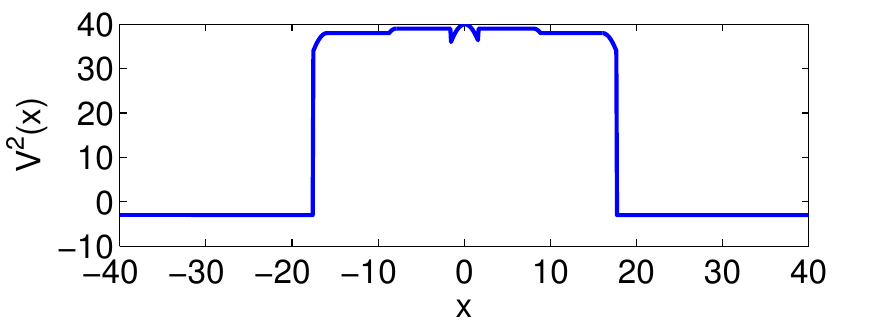}

\caption{\footnotesize Optimal sum of rewards (value) 
$V^t(x)$ for $b = 0 \, 
(\false)$ for time horizons (i.e., decision stages remaining) $t=0$,
$t=1$, and $t=2$ on the \textsc{Continuous Action}  \MarsRover\ problem.  For $x \in [-2,2]$, the
rover automatically takes a picture and receives a reward quadratic in
$x$.  We initialized $V^0(x,b) = R(x,b)$; for $V^1(x)$, the rover achieves
non-zero value up to $x = \pm 12$ and for 
$V^2(x)$, up to $x = \pm 22$.}
\label{fig:opt_graph}

\end{minipage}

\begin{minipage}[b]{0.8\linewidth}
\includegraphics[width=1\textwidth]{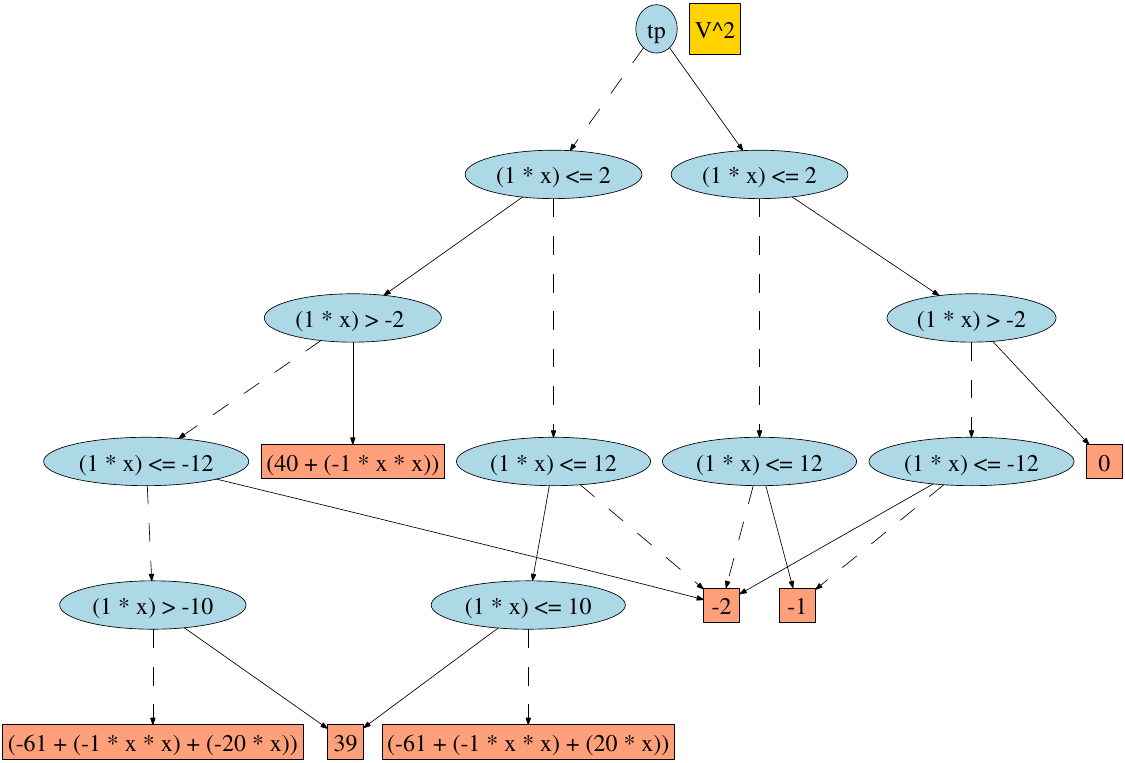}
\vspace{6mm}

\caption{\footnotesize Optimal value function $V^2(x)$ for the
\textsc{Continuous Action}  \MarsRover\ problem represented as an (XADD) equal to the second diagram on the left. 
To evaluate
$V^2(x)$ for any state $x$, one simply traverses the diagram in a
decision-tree like fashion until a leaf is reached where the
non-parenthetical expression provides the \emph{optimal value} and the
parenthetical expression provides the \emph{optimal policy} 
($y = \pi^{*,2}(x)$) to achieve value $V^2(x)$.}
\label{fig:opt_val_pol}
\vspace{-5mm}
\end{minipage}

\end{figure}
Assuming our goal is to amplify the long haul \emph{value} $V$ (i.e.,
the amount of remunerations got over a limitless skyline of activities), then
we can compose the ideal worth feasible from a given state in \MarsRover\
as an element of state factors.

The worth capability is piecewise and non-straight, and it contains non-rectangular choice
limits like $4 - x^2 - y^2\geq 0$. Figure~\ref{fig:opt_graph} presents the 0-, 1-, and 2-step time skyline answers for this issue; further, in representative structure, we show both the 1-step time skyline esteem capability and relating ideal arrangement in Figure~\ref{fig:opt_val_pol}. Here, the piecewise idea of the progress and award capability prompts piece-shrewd construction in the worth capability and strategy. However notwithstanding the instinctive and straightforward nature of this outcome, we know nothing about earlier strategies that can create such definite arrangements.

\paragraph{\WaterReservoir}
Supply the board is all around concentrated in
the OR writing \cite{Mahootchi2009,Yeh1985}. The key constant choice is the ticket
much slipped by time $e$ to
\emph{drain} (or \emph{not drain}) every supply to boost
power income over the choice stage skyline while keeping away from
supply flood and sub-current. Give a role as a CA-HMDP, we
accept SVI gives the principal approach equipped for determining
an accurate shut structure non-nearsighted ideal approach
for all levels.

We look at a 2-repository issue with
separate levels $(l_1,l_2)\in [0,\infty]^2$ with remuneration punishments for
flood and sub-current and a prize increase direct in the passed time $e$ for
power created in periods when the $\mathit{drain}(e)$ activity
channels water from $l_2$ to $l_1$ (the other activity is
$\mathit{no}$-$\mathit{drain}(e)$); we expect deterministic precipitation
recharging and present the prize capability as:
\vspace{-4mm}
\begin{figure*}[h]
\centering
\footnotesize
\begin{align}
R & = 
\begin{cases}
((50 - 200e) \leq l_1 \leq (4500 - 200e)) \wedge ((50 + 100e) \leq l_2 \leq (4500 + 100e)) & : e \\
((50 + 300e) \leq l_1 \leq (4500 + 300e)) \wedge ((50 - 400e) \leq l_2 \leq (4500 - 400e)) & : 0 \\
\text{otherwise} & : -\infty \\
\end{cases} 
\end{align}
\end{figure*}
The transition function for levels of the $\mathit{drain}$ action is defined below. Note that for the $\mathit{no}$-$\mathit{drain}$ action, the $\mathit{500 * e}$ term is not involved.
{
\begin{align*}
l_1' & =(400 * e + l_1 -700 * e + 500 * e) \\
l_2'& =(400 * e + l_2 - 500 * e) \\
\end{align*}}

Similar to the discrete version of the \InventoryControl problem in the introduction, the DA-HMDP setting for these two problems defines discrete actions by partitioning the action space of each domain into $i$ number of slices. For example the \MarsRover\ problem with 2 actions which are the lower and upper bounds on $a$ ($a_1=-10, a_2 = 10$) and the transition and reward functions are defined according to these constant values.
We now provide the empirical results obtained from implementing our algorithms.

\subsection{Results}

For both the \textsc{Discrete Action} and \textsc{Continuous Action} \MarsRover\ domains, 
we have run experiments to evaluate our SDP solution 
in terms of time and space cost while varying the horizon and problem size.

\begin{figure*}[t]
\centering

\includegraphics[width=0.45\textwidth]{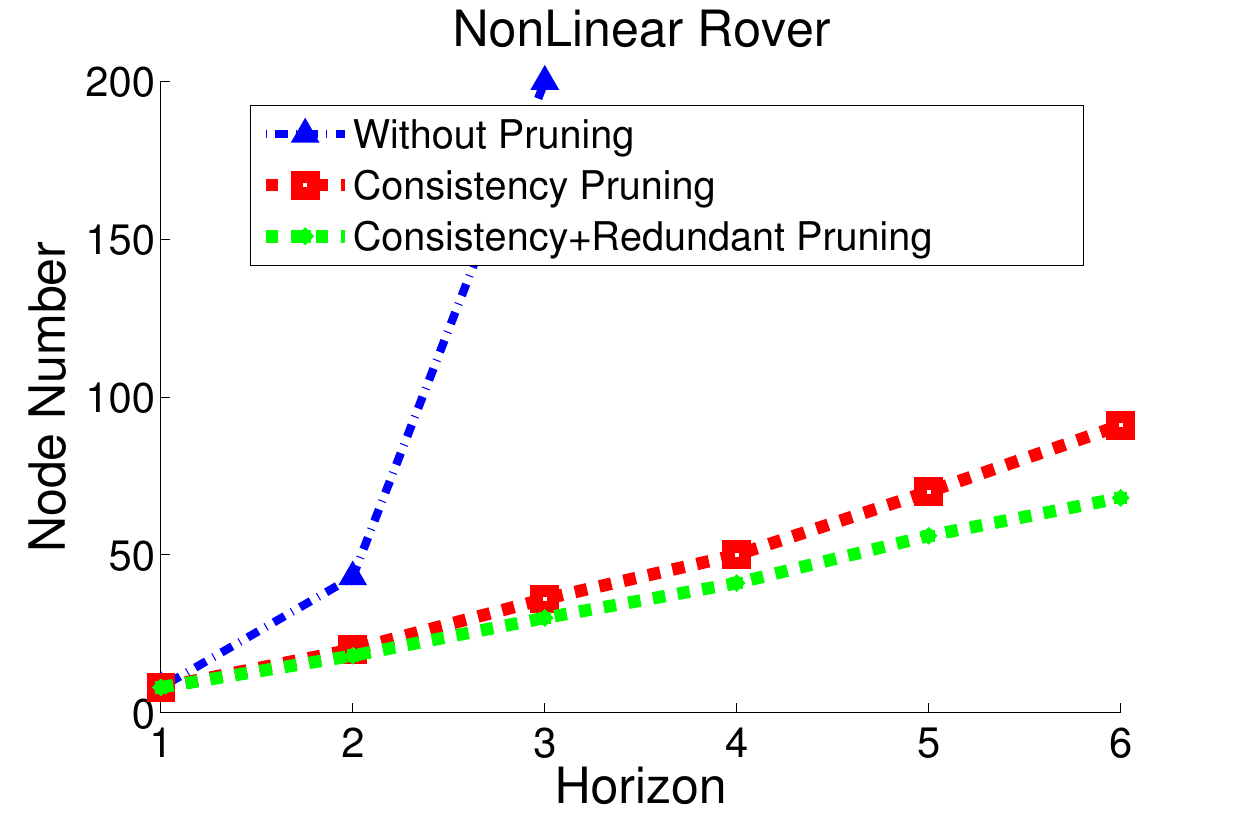}

\includegraphics[width=0.45\textwidth]{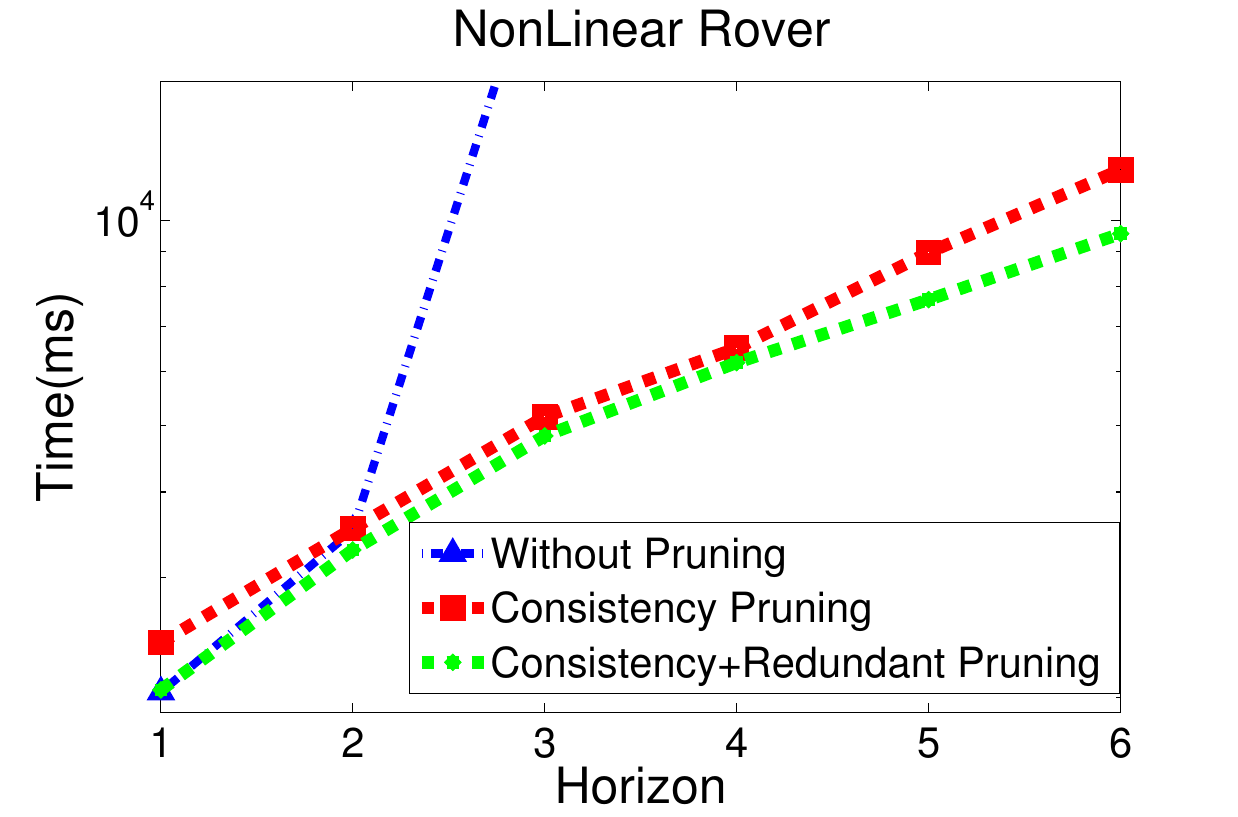}
\vspace{-3mm}
\caption{
Space (\# XADD nodes in value function) and
time for different iterations (horizons) of SDP on Nonlinear \textsc{Continuous Action}  \MarsRover\ with 4 different results based on pruning techniques. Results are shown for  the XADD 
 with no pruning technique, with only consistency checking (using LP-solver) and with both the consistency and redundancy checking (using LP-Solver and SAT-Solver) with a numerical precision heuristic.} 
\label{fig:roverTS}

\end{figure*}

\begin{figure}[h]
\vspace{2mm}
\centering
\includegraphics[width=0.45\textwidth]{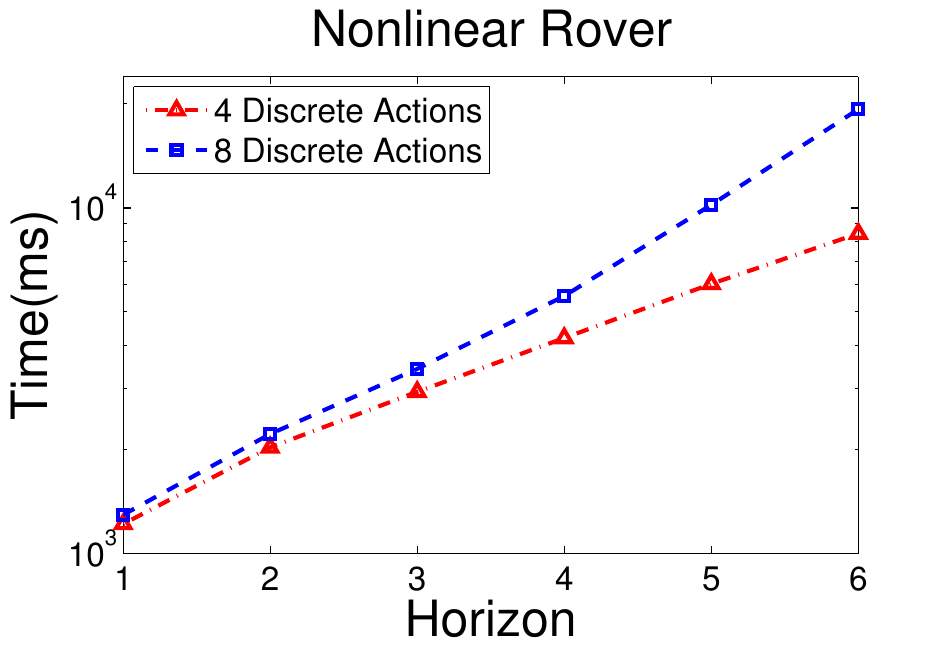}
\hspace{2mm}
\includegraphics[width=0.45\textwidth]{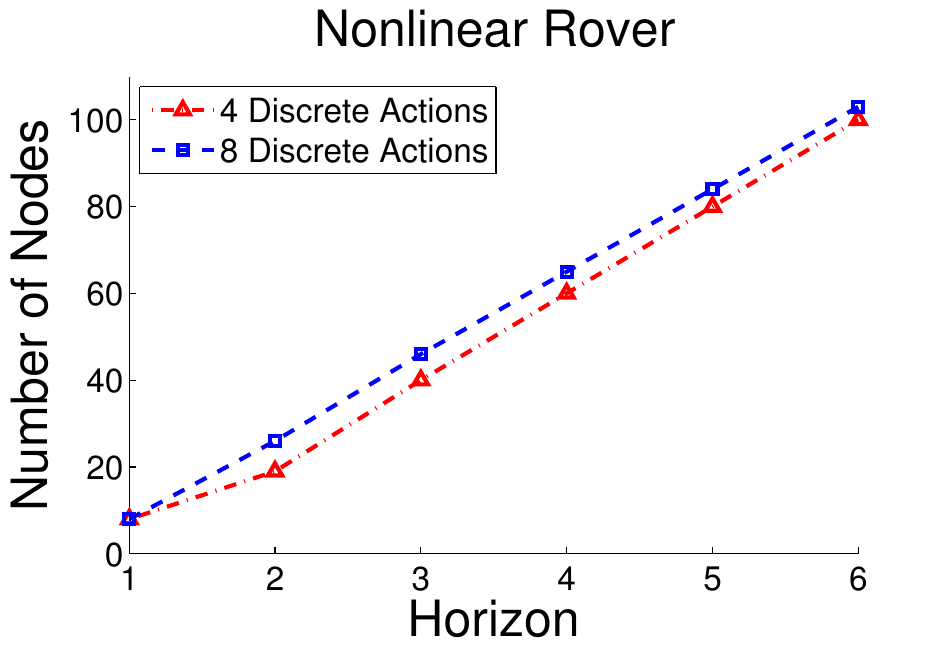}
\vspace{-2mm}
\caption{
Space and elapsed time vs. horizon for a fixed discritization of 4 and 8 for the \MarsRover to show how the number of nodes and time increases for each horizon. 
}
\label{fig:roverDisTS}
\vspace{-5mm}
\end{figure}

\begin{figure}[h]
\vspace{-2mm}
\centering
\includegraphics[width=0.45\textwidth]{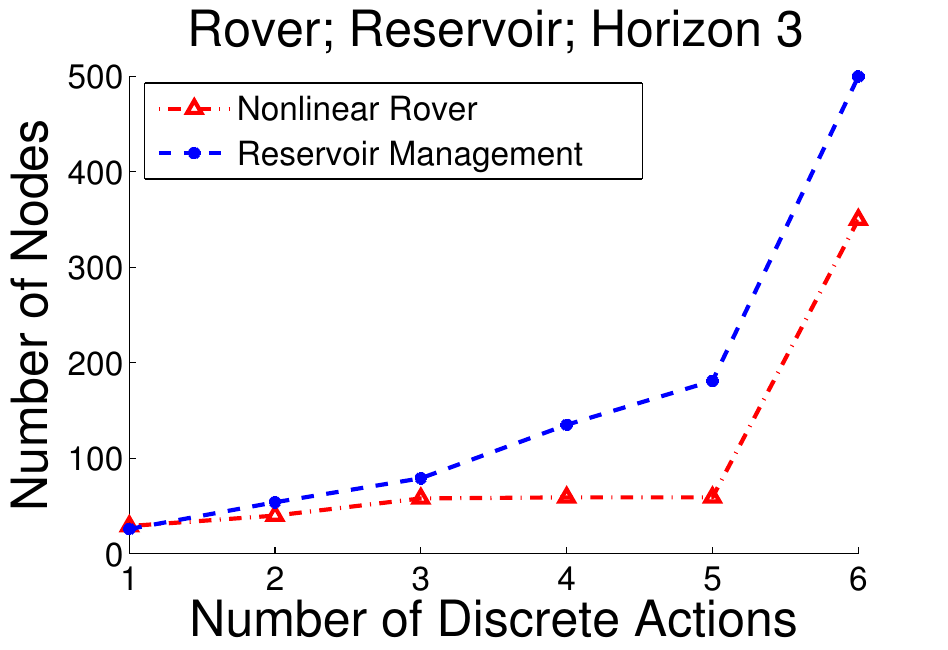}
\hspace{2mm}
\includegraphics[width=0.45\textwidth]{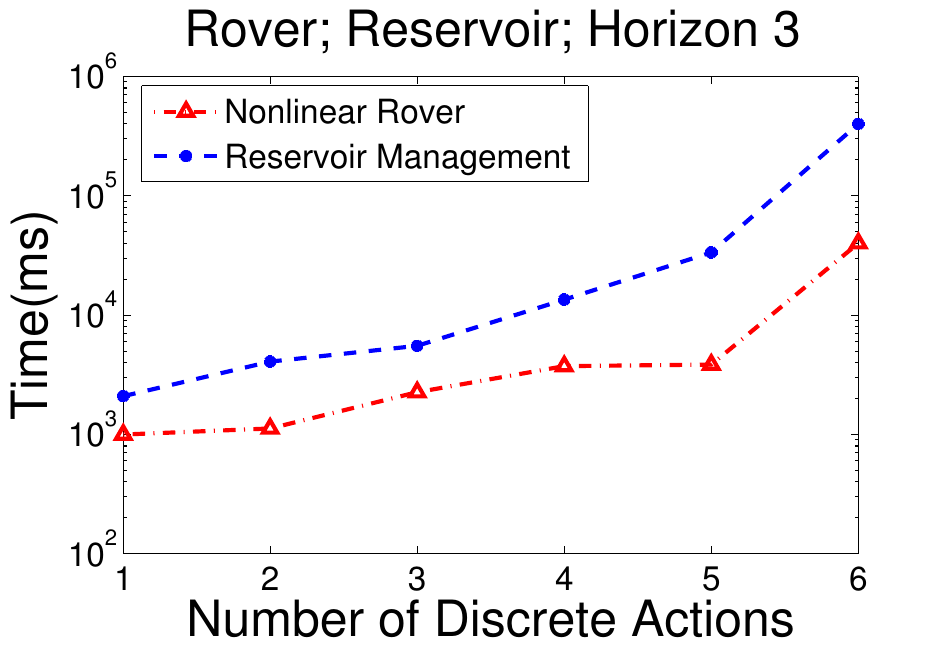}
\vspace{-2mm}
\caption{
Space and elapsed time vs. horizon for a fixed horizon of 3 and different problem sizes for the \textsc{Discrete Action} \MarsRover\ and the \WaterReservoir. 
}
\label{fig:roverDisSize}
\vspace{-5mm}
\end{figure}

\begin{figure}[h]
\vspace{2mm}
\centering
\includegraphics[width=0.45\textwidth]{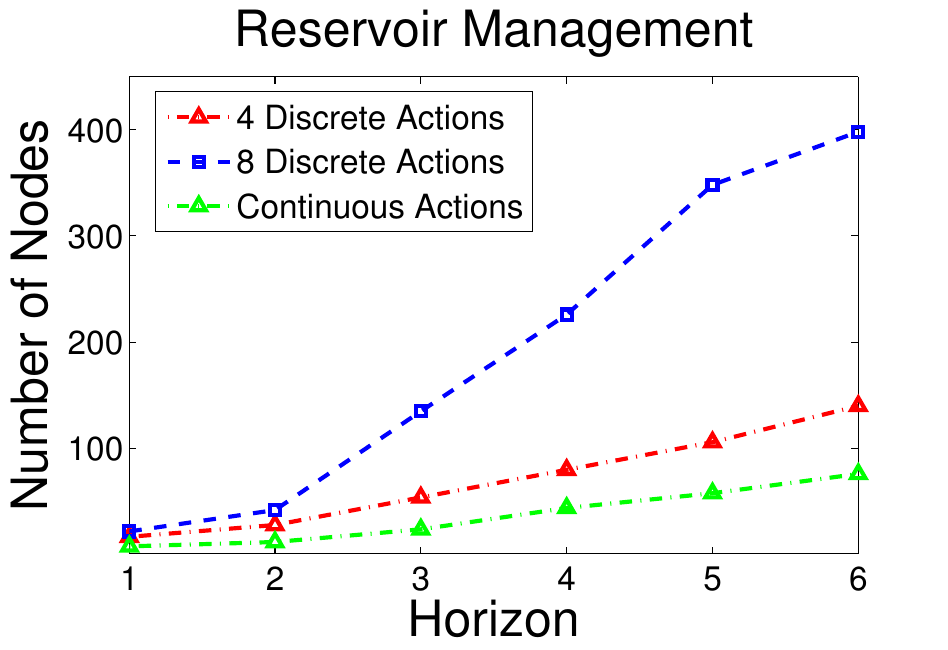}
\hspace{2mm}
\includegraphics[width=0.45\textwidth]{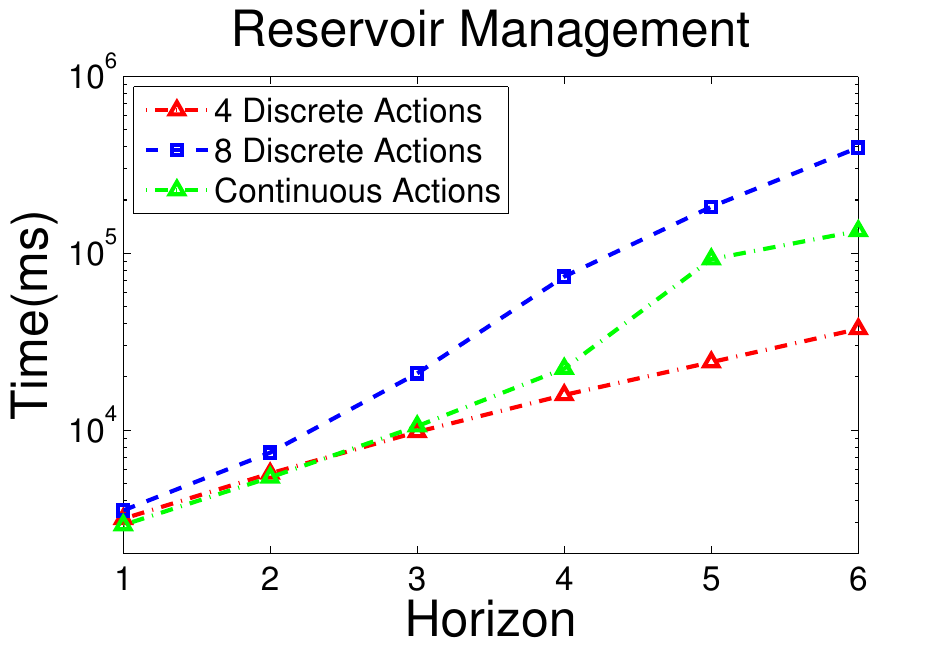}
\vspace{-2mm}
\caption{
Space and elapsed time vs. horizon for a fixed discritization of 4 and 8 for the discrete action \WaterReservoir. To compare the time and space of continuous action \WaterReservoir is given. }
\label{fig:resDisTS}
\vspace{-5mm}
\end{figure}

\begin{figure*}[h]
\vspace{-2mm}
\centering
\includegraphics[width=0.29\textwidth]{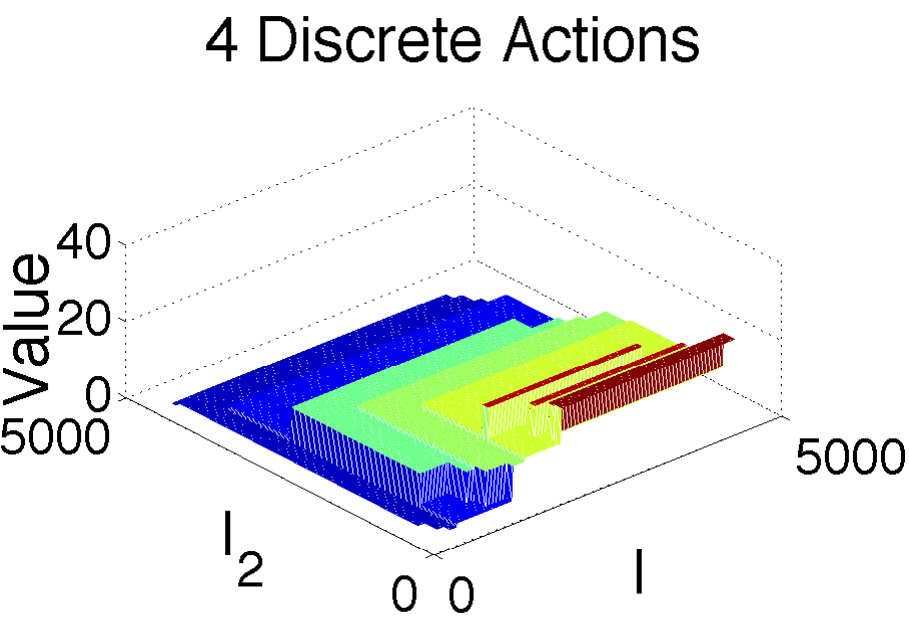}
\includegraphics[width=0.29\textwidth]{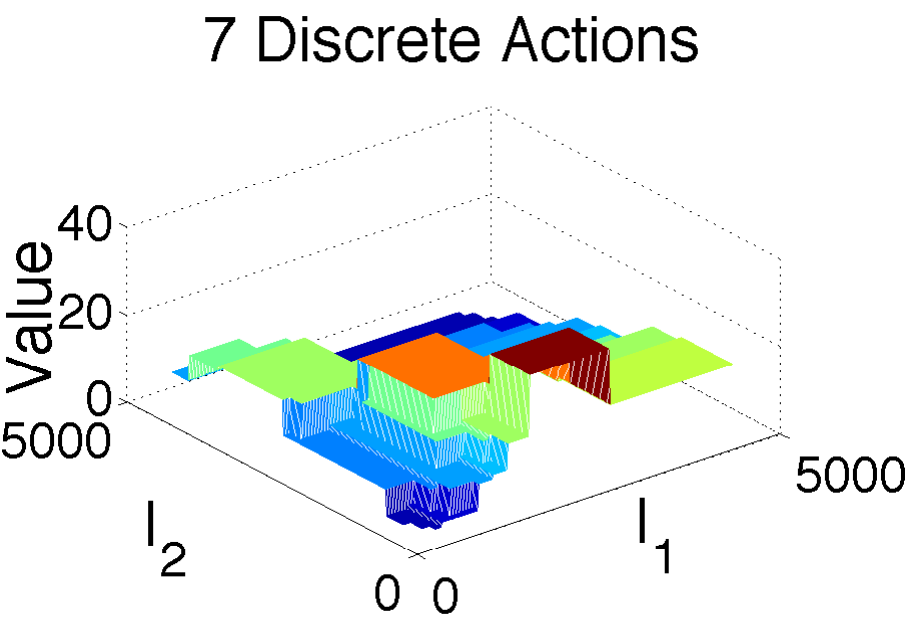}
\includegraphics[width=0.29\textwidth]{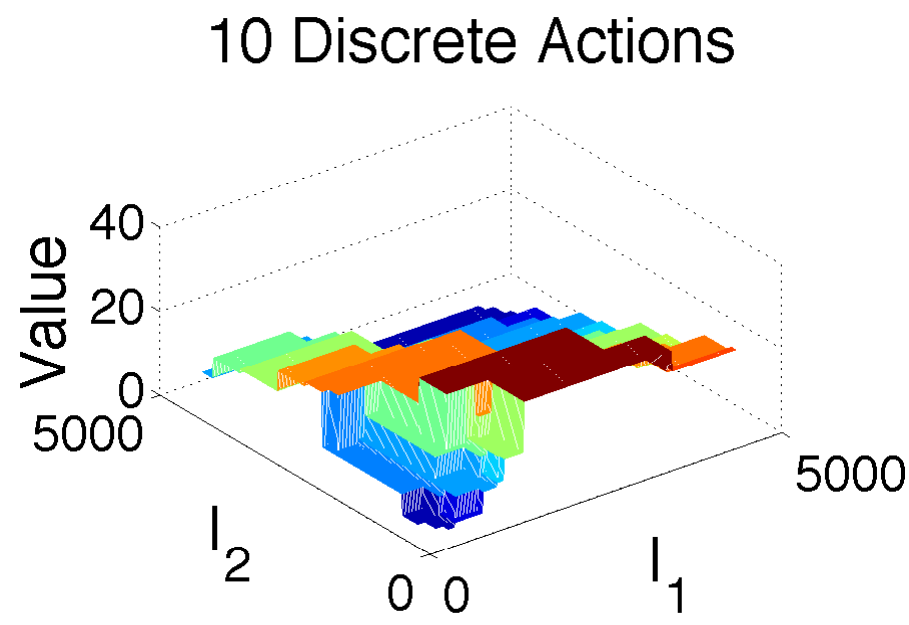}
\vspace{-3mm}
\caption{
4, 7 and 10 \textsc{Discrete Actions} of the \WaterReservoir problem for different water levels in iteration $V^3$. Each discretization draws the value closer to that of the continuous action \WaterReservoir in Figure~\ref{fig:v2plots}.  
}
\label{fig:discreteplots}
\vspace{-4mm}
\end{figure*}
For the \textsc{Continuous Action} \MarsRover problem, we present the time and space analysis for the problem description defined in the introduction in Figure~\ref{fig:roverTS}. Here three evaluations are performed based on the pruning algorithms in the previous section. We note that without the LP-Solver Algorithm~\ref{algPrune}, the algorithm can not go beyond the forth iteration as it produces many inconsistent nodes. Next we use the LP-solver for consistency pruning and then use redundancy checking. Note that this approach has an added numerical heuristic that omits similar branches. With the full pruning of Algorithm~\ref{algRedundant} less than 10 \% node reduction is achieved but due to the calls made to the SAT-Solver the time increases to more than 10 times of that of the heuristic pruning.  Also note that even the heuristic approach has very little node reduction can be gained from redundancy pruning. This suggests that in some domains using redundancy pruning in not efficient and should be omitted from the final results. 
We will only show results for the  \WaterReservoir and \InventoryControl problem with consistency pruning. 

Next we present the analysis of the \textsc{Discrete Action} \MarsRover\ domains. Figure~\ref{fig:roverDisTS} shows how time and space costs of different horizons increases for a fixed action discritization of 4 and 8.
Note that one of the caveats of using a discrete setting is defining the actual discrete actions. While in the continuous setting an action is defined between a large high and low range (e.g $a \in$ [$-1000000,1000000$]) allowing the SDP algorithm to choose the best possible action among all the answers, for a discrete setting, knowing the range to discretize the action becomes very important. As an example in the rover description, allowing actions to be far from the center (e.g. $a \in$ [$-20,20$]) does not result in a converged solution. Finding this range is one of the drawbacks of using a discrete setting. 

On the other hand, the level of discretization within this range also is very important. To show this visually Figure shows the results of the third iteration for 6 different discretizations compared to the continuous result. This figure proves the need to partition the action space within the predefined range more finely, but as Figure~\ref{fig:roverDisSize} shows  the effect of varying the  problem size for the fixed horizon of 3, leads to increasing time and space costs. Here results of different level of discretization is similar for \WaterReservoir.

Also Figure~\ref{fig:resDisTS} presents the number of nodes and time for different horizons of 4 and 8 discrete actions compared to the continuous actions in the \WaterReservoir domain. While the first two iterations the time and nodes follows closely as the number of actions doubles, for higher iterations more time and space is required. The number of nodes are lower for the continuous compared to both discretizations but the time elapsed is higher than the 4-discrete actions due to the complexity of the continuous action maximization.

\begin{figure*}[h]
\vspace{-2mm}
\centering
\includegraphics[width=0.29\textwidth]{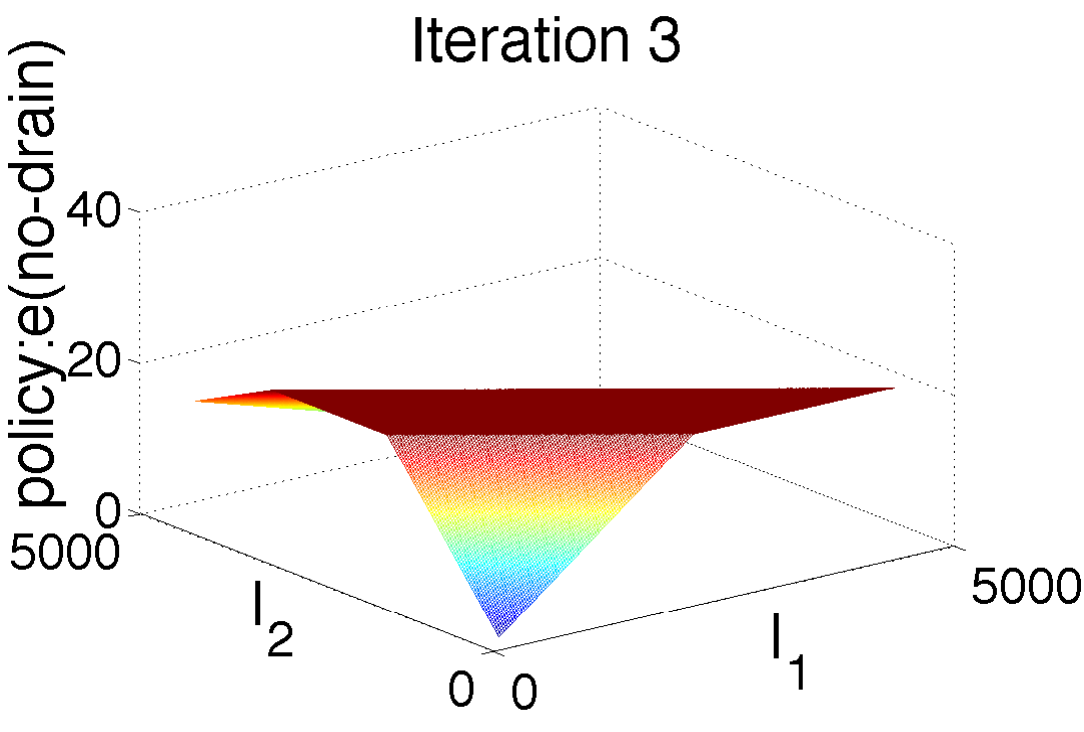}
\includegraphics[width=0.29\textwidth]{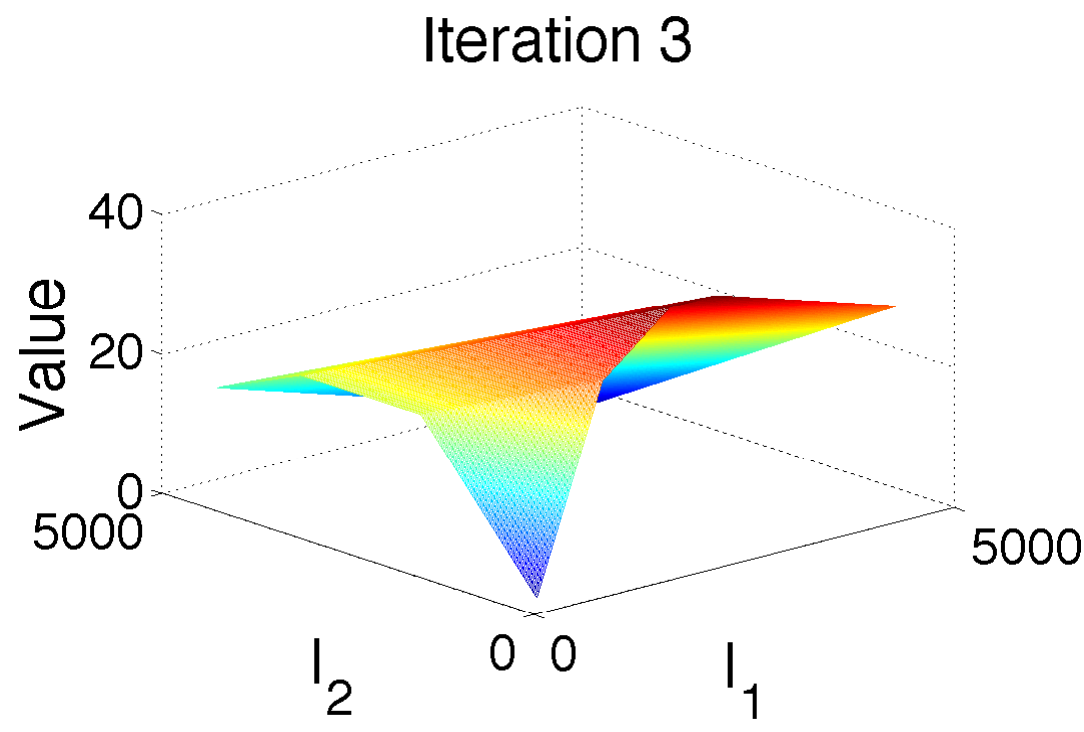}
\includegraphics[width=0.29\textwidth]{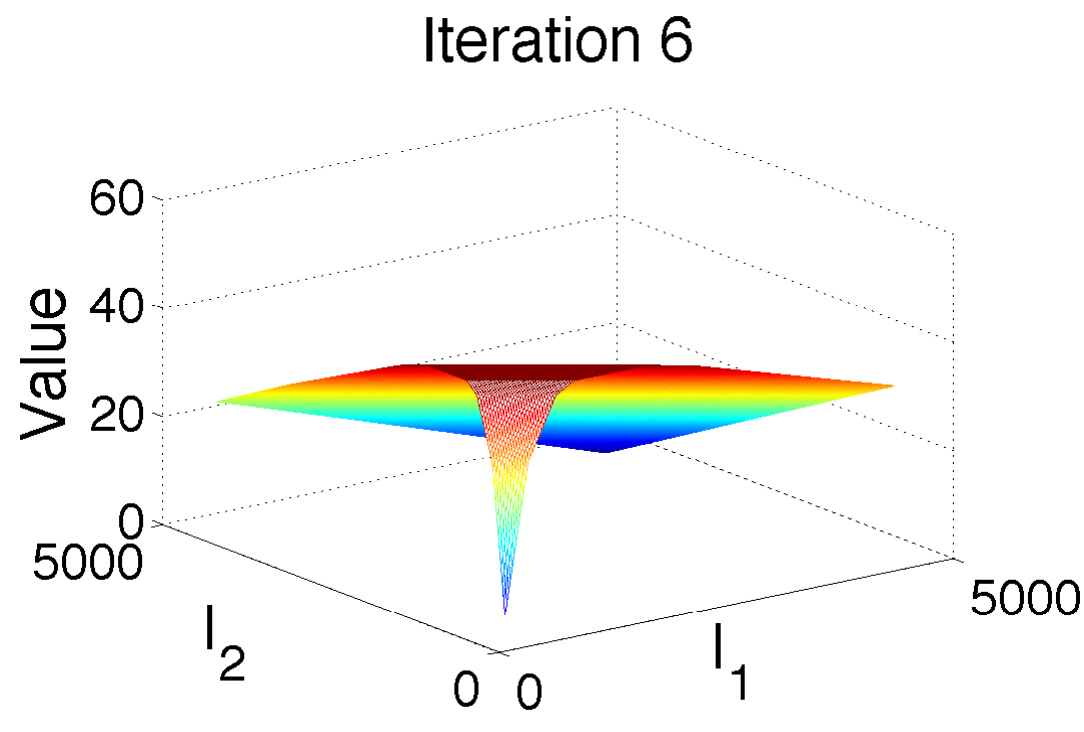}
\vspace{-3mm}
\caption{
\WaterReservoir: 
{\it (left)} Policy $\mathit{no}$-$\mathit{drain}(e)=\pi^{3,*}(l_1,l_2)$ 
showing on the z-axis the elapsed time $e$ that should be execu                                          ted 
for $\mathit{no}$-$\mathit{drain}$ conditioned on the states; 
{\it (middle)} $V^3(l_1,l_2)$; 
{\it (right)} $V^6(l_1,l_2)$.
}
\label{fig:v2plots}
\end{figure*}

Figure~\ref{fig:discreteplots} shows 3 discretizations of the \WaterReservoir problem. For the left figure 4 discrete actions, the middle figure 7 actions and the right figure uses 10 discrete actions to represent the value of the third iteration according to water levels $l_1$ and $l_2$. The figures suggest that finer grain discretization results in better results, closer to that of the continuous action value in middle picture of Figure~\ref{fig:v2plots}. 

In Figure~\ref{fig:v2plots}, we show a plot of
the ideal shut structure strategy at $h=3$: the arrangement interleaves $\mathit{drain}(e)$ and $\mathit{no}$-$\mathit{drain}(e)$ where even skylines are the last option;
here we see that we try not to deplete for the longest passed time $e$
when $l_2$ is low (trust that downpour will recharge) and $l_1$ is high (depleting
water into it could spill over it). $V^3(l_1,l_2)$ and $V^6(l_1,l_2)$
show the movement of union from skyline $h=3$ to $h=6$ - - -
low degrees of $l_1$ and $l_2$ permit the framework to create power
for the longest absolute slipped by time north of 6 choice stages.

\begin{figure}[h]
\centering
\includegraphics[width=0.45\textwidth]{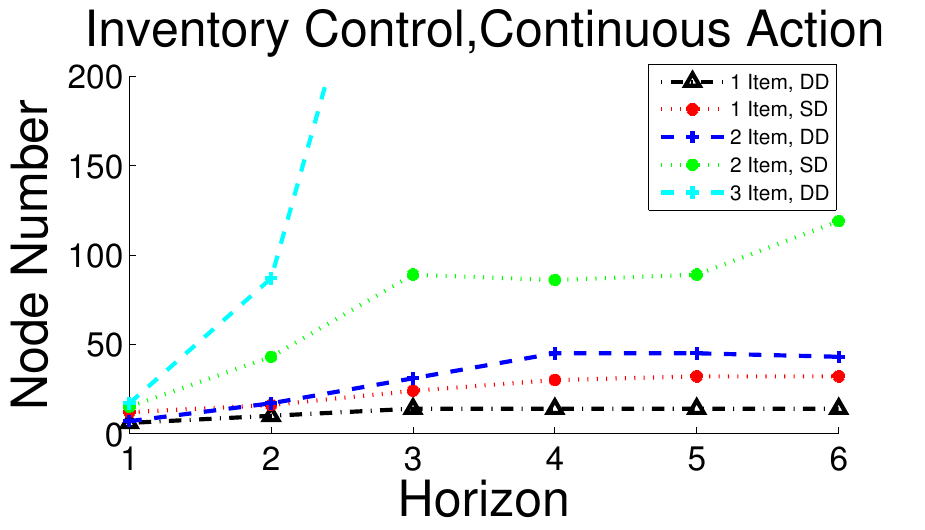}
\hspace{2mm}
\includegraphics[width=0.45\textwidth]{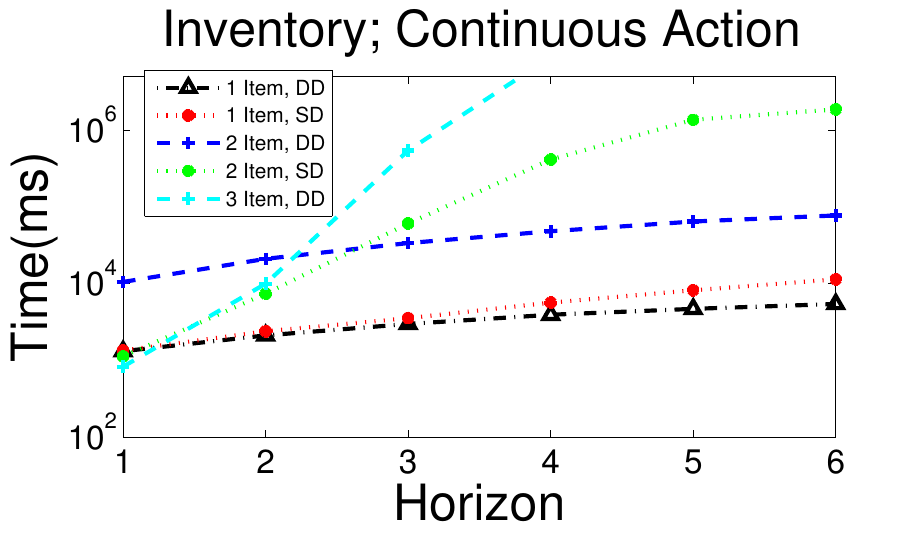}
\caption{\textsc{Continuous Action} \InventoryControl: Space and time vs. horizon.}
\label{fig:invC}
\end{figure}

In Figure~\ref{fig:invC}, we give a reality examination of
deterministic-and stochastic-request (resp. DD and SD) variations of the
SCIC and MJCIC issue for up to three things (a similar size of
issues frequently concentrated in the OR writing); for each number of things
$n \in \{ 1,2,3 \}$ the state (stock levels) is $\vec{x} \in
[0,\infty]^n$ and the activity (reorder sums) is $\vec{y} \in
[0,\infty]^n$. Orders are made at one month stretches and we settle
for a skyline up to $h=6$ months.

\begin{figure}[h]
\centering
\includegraphics[width=0.45\textwidth]{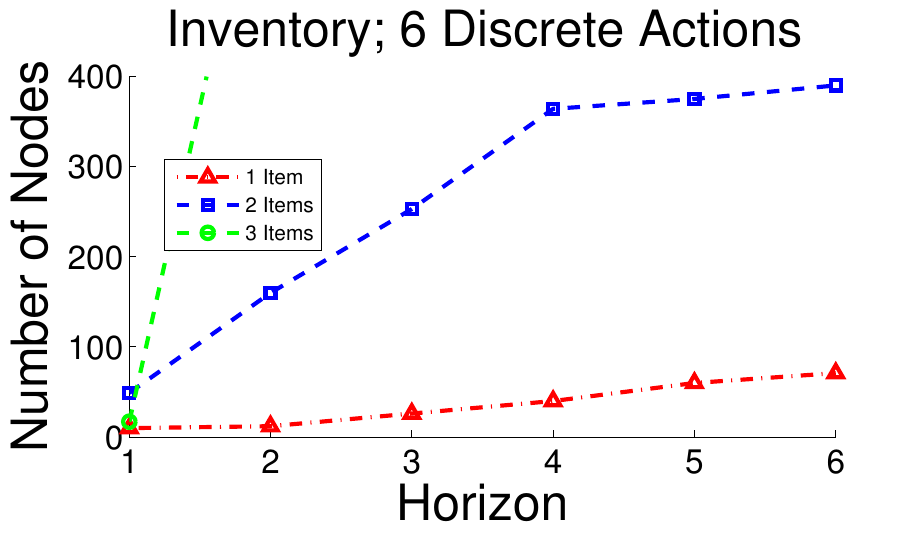}
\hspace{2mm}
\includegraphics[width=0.45\textwidth]{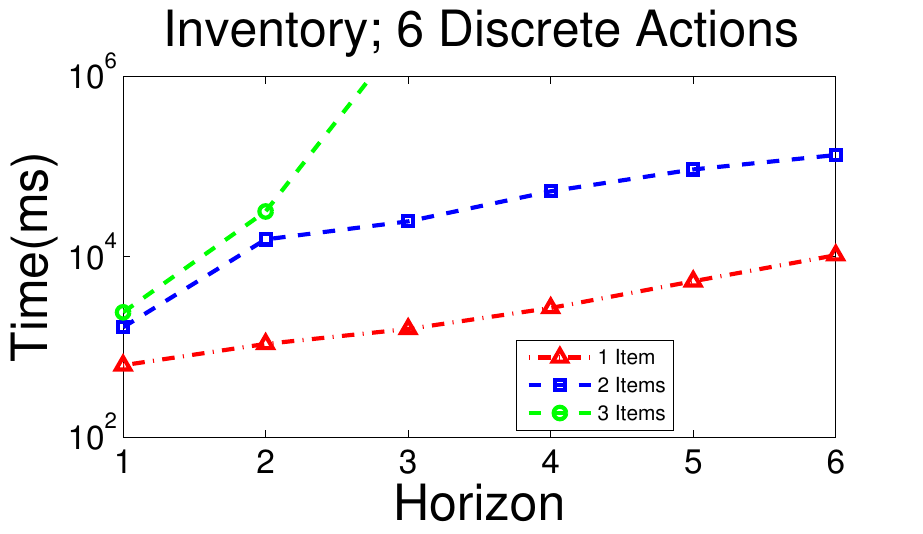}
\caption{\textsc{Discrete Action} \InventoryControl: Space and time vs. horizon for 6 discrete actions. Results are shown for various continuous items in the inventory.}
\label{fig:invD6}
\end{figure}

While addressing for bigger quantities of
things and SD (as opposed to DD) both increment existence,
the arrangements rapidly arrive at calmness demonstrating underlying
union.
Figure~\ref{fig:invD6} addresses the existence for the deterministic \textsc{Discrete Action} \InventoryControl for a discretization of 6 activities and different stock things for up to $h=6$ skylines. While the quantity of things influences both existence, in any event, for 6 discrete activities the 3-thing stock will have remarkable reality for the second skyline onwards. The explanation is behind the high aspects, for a 3-thing stock of 6 discrete activities, a sum of $6 \times 6 \times 6$ activities are required!

\begin{figure}[h]
\centering
\includegraphics[width=0.45\linewidth]{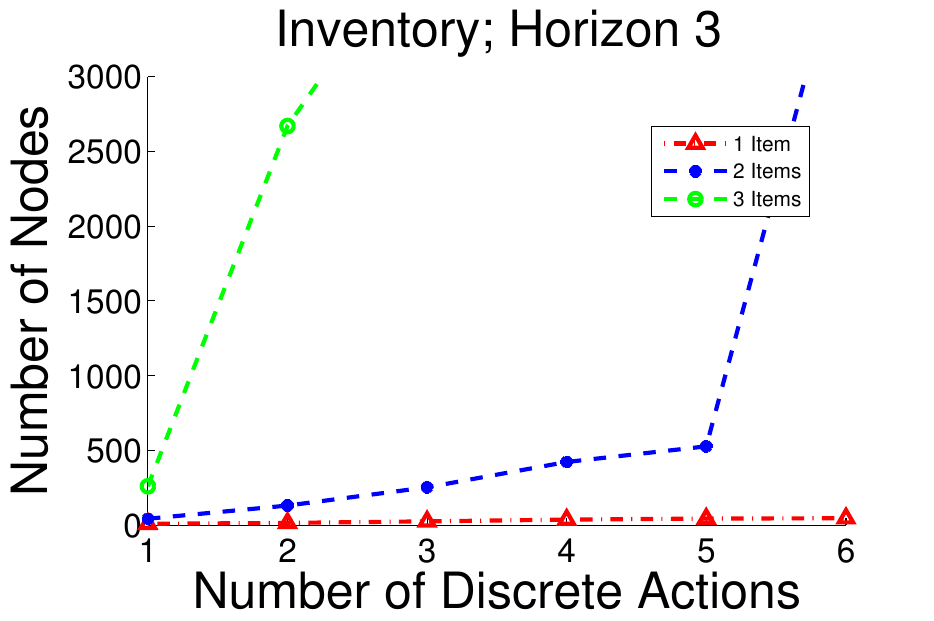}
\hspace{2mm}
\includegraphics[width=0.45\linewidth]{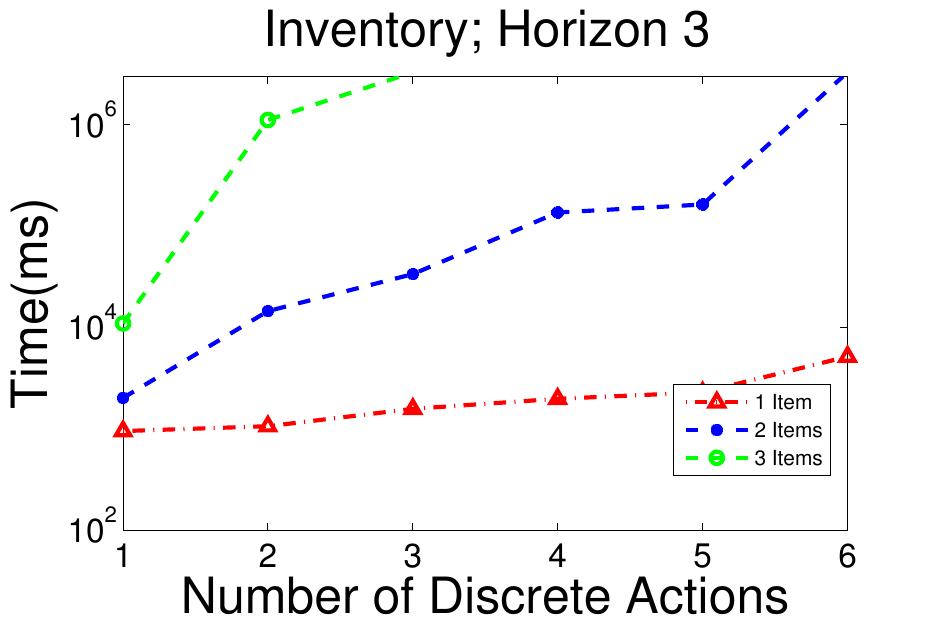}
\caption{Space and elapsed time vs. horizon for a fixed horizon of 3 and different problem sizes for the \textsc{Discrete Action} \InventoryControl.}
\label{fig:invH3}
\end{figure}

To show the effect of different number of action discretizations for the \InventoryControl Figure~\ref{fig:invH3} is used to analyse the time and space for a fixed horizon of $h=3$ and inventory items of $1,2$ and 3. As the number of discrete actions increases, the time and space grows almost linearly for the 1 and 2 item inventory but again for a 3-item inventory this causes major space and time limitations. 
\section{Related Work}

The most pertinent vein of related work for DA-HMDPs is that of \cite{feng04} and \cite{li05} which can perform careful unique programming on
HMDPs with rectangular piecewise straight prize and change capabilities
that are delta capabilities. While SDP can tackle these equivalent issues,
it eliminates both the rectangularity and piecewise limitations on the
award and worth capabilities, while holding precision.
Heuristic inquiry approaches with formal assurances
like HAO* \cite{hao09} are an appealing future augmentation of SDP;
HAO* right now utilizes the strategy for \cite{feng04}, which could, as a matter of fact
be straightforwardly supplanted with SDP. While \cite{penberthy94} has thought of
general piecewise capabilities with straight limits (and truth be told,
we get our direct pruning come closer from this paper), this work
simply applied to completely deterministic settings, not HMDPs.

Other work has examined restricted HMDPS having only one ceaseless
state variable. Obviously rectangular limitations are insignificant with
just a single consistent variable, so it isn't is business as usual that more
progress has been made in this limited setting. One persistent
variable can be helpful for ideal answers for time-subordinate MDPs
(TMDPs) \cite{boyan01}. Or on the other hand stage changes can be utilized to
with no obvious end goal in mind surmised one-layered persistent disseminations
prompting a limited estimation approach for inconsistent single persistent
variable HMDPs \cite{phase07}.
While this work can't deal with inconsistent stochastic
commotion in its persistent conveyance, it really does precisely tackle HMDPs
with different constant state aspects.

There are various general HMDP estimation
approaches that utilization surmised direct programming \cite{kveton06}
or on the other hand examining in a support learning style approach \cite{munos02}.
As a rule, while estimate strategies are very encouraging in
practice for HMDPS, the goal of this paper was to push
the limits of \emph{exact} arrangements; notwithstanding, in some sense,
we accept that more expressive precise arrangements may likewise illuminate
better approximations, e.g., by permitting the utilization of information structures
with non-rectangular piecewise parts that permit higher constancy
approximations.

Concerning CA-HMDPs, there has been earlier work in charge hypothesis. The field of straight quadratic Gaussian (LQG) control \cite{lqgc} which utilize direct elements with consistent activities, Gaussian clamor, and quadratic
reward is generally firmly related. Nonetheless, these precise arrangements do
not stretch out to discrete and consistent frameworks with \emph{piecewise}
elements or prize.

Joining this work with starting state centered methods \cite{hao09}
what's more, engaged approximations that exploit ideal worth
structure \cite{apricodd} or further
away from home \cite{munos02,kveton06,phase07} are promising headings for
future work.

\section{Concluding Remarks}

In this paper, we acquainted another emblematic methodology with tackling nonstop issues in HMDPs precisely. On account of discrete activities and nonstop states, utilizing erratic
reward capabilities and expressive nonlinear progress capabilities far surpasses the specific arrangements conceivable with existing HMDP
solvers.
Concerning ceaseless states and activities, a key commitment is that of \emph{symbolic compelled
optimization} to take care of the persistent activity amplification issue. We
accept this is the primary work to propose ideal shut structure
answers for MDPs with \emph{multivariate} nonstop state \emph{and}
activities, discrete commotion, \emph{piecewise} direct elements, and
\emph{piecewise} straight (or limited \emph{piecewise} quadratic)
reward; further, we accept our exploratory outcomes are the first
careful answers for these issues to give a shut structure ideal
strategy for every single (nonstop) state.
While our strategy isn't versatile for 100's of things, it actually addresses
the primary general accurate arrangement techniques for capacitated multi-stock control issues.
What's more, albeit a straight or quadratic prize is very restricted yet it has seemed valuable for single consistent asset or constant time issues like the water repository issue.

With an end goal to make SDP useful, we additionally presented
the original XADD information structure for addressing erratic piecewise
emblematic worth capabilities and we tended to the difficulties that
SDP prompts for XADDs, like the requirement for reordering and pruning the choice
hubs after certain activities. These are significant commitments
that have added to another degree of expressiveness for HMDPS
that can be precisely addressed.

There are various roads for future examination. Most importantly, it is
significant analyze what speculations of the change capability utilized
in this work would in any case allow shut structure accurate arrangements. In wording
of better versatility, one road would investigate the utilization of introductory
state centered heuristic hunt based esteem cycle like
HAO* \cite{hao09} that can be promptly adjusted to utilize SDP. Another
road of exploration is adjust the languid estimate approach
of \cite{li05} to rough HMDP esteem capabilities as piecewise
direct XADDs with straight limits that might take into account better
approximations than current portrayals that depend on rectangular
piecewise capabilities. Similarly, thoughts from
APRICODD \cite{apricodd} for limited guess of discrete ADD
esteem capabilities by combining leaves could be summed up to XADDs.
By and large the advances made by this stir open up various
potential novel exploration ways that we accept may help make
fast advancement in the field of choice hypothetical preparation
with discrete and ceaseless state.

With the ongoing answer for ceaseless states and activities, we can apply our techniques to certifiable information from the Stock writing with additional careful advances and rewards. Completely stochastic circulations are expected for these issues which is a significant future heading by in-participating a commotion boundary in the models.
Likewise we have investigated esteem emphasis for the two issues, settling the emblematic approach cycle calculation for issues with straightforward strategies can end up being successful in specific spaces.
The other promising course is to expand the ongoing definite answer for non-direct capabilities and addressing polynomial conditions utilizing computational mathematical strategies.

\bibliographystyle{IEEEtran}
\bibliography{exactsdp}

\end{document}